\documentclass[11pt]{amsart}
\usepackage{t1enc}
\usepackage[latin1]{inputenc}
\usepackage[english]{babel}
\usepackage{amsmath,amsthm}
\usepackage{amsfonts}
\usepackage{latexsym}
\usepackage[dvips]{graphicx}
\usepackage{graphicx}
\usepackage{float}
\usepackage[natural]{xcolor}
\usepackage{algorithm}
\usepackage{algorithmic}
\usepackage{enumerate}
\usepackage{appendix}
\usepackage{multirow}
\usepackage{xcolor}
\usepackage[colorlinks,linkcolor=blue]{hyperref}
\usepackage{lineno}
\usepackage{setspace}

\usepackage{fullpage}

\usepackage{listings}
\newtheorem{theorem}{Theorem}\numberwithin{theorem}{section}
\newtheorem{lemma}[theorem]{Lemma}\numberwithin{theorem}{section}

\numberwithin{theorem}{section}
\newtheorem{fact}{Fact}

\newtheorem{definition}[theorem]{Definition}

\newtheorem{conjecture}[theorem]{Conjecture}
\newtheorem{claim}[theorem]{Claim}

\def\eps{\varepsilon}

\def\int{\textrm{int}}

\begin{document}

\onehalfspacing
\title{Rainbow spanning structures in graph and hypergraph systems}

\author{Yangyang Cheng}
\address{YC, BW and GW. School of Mathematics, Shandong University, Jinan, China\\
Email: \texttt{(YC) mathsoul@mail.sdu.edu.cn, (BW) binwang@mail.sdu.edu.cn, (GW) ghwang@sdu.edu.cn}}
\author{Jie Han}
\address{JH. School of Mathematics and Statistics, Beijing Institute of Technology, Beijing,China, Email: \texttt{han.jie@bit.edu.cn}.}
\author{Bin Wang}
\author{Guanghui Wang}


\begin{abstract}
We study the following rainbow version of subgraph containment problems in a family of (hyper)graphs, which generalizes the classical subgraph containment problems in a single host graph.
For a collection $\emph{\textbf{G}}=\{G_1, G_2,\ldots, G_{m}\}$ of not necessarily distinct $k$-graphs on the same vertex set $[n]$, a (sub)graph $H$ on $[n]$ is rainbow if there exists an injection $\varphi: E(H)\rightarrow[m]$ such that $e\in E(G_{\varphi(e)})$ for each $e\in E(H)$.
Note that if $|E(H)|=m$, then $\varphi$ is a bijection and thus $H$ contains exactly one edge from each $G_i$.

Our main results focus on rainbow clique-factors in (hyper)graph systems with minimum $d$-degree conditions.
Specifically, we establish the following:
\begin{enumerate}
\item A rainbow analogue of an asymptotical version of the Hajnal--Szemer\'{e}di theorem, namely, if $t\mid n$ and $\delta(G_i)\geq(1-\frac{1}{t}+\varepsilon)n$ for each $i\in[\frac{n}{t}\binom{t}{2}]$, then $\emph{\textbf{G}}$ contains a rainbow $K_t$-factor;
\item Essentially a minimum $d$-degree condition forcing a perfect matching in a $k$-graph also forces rainbow perfect matchings in $k$-graph systems for $d\in[k-1]$.
\end{enumerate}
The degree assumptions in both results are asymptotically best possible (although the minimum $d$-degree condition forcing a perfect matching in a $k$-graph is in general unknown).
For (1) we also discuss two directed versions and a multipartite version.
Finally, to establish these results, we in fact provide a general framework to attack this type of problems, which reduces it to subproblems with \emph{finitely many} colors.

\bigskip


\end{abstract}

\maketitle

\section{Introduction}

\subsection{Rainbow extremal graph theory}
A natural variant of the extremal problems concerns rainbow substructures in edge-colored graphs.
From our knowledge two types of host graphs have been studied: one class is the properly edge-colored graphs, which was first considered by Keevash, Mubayi, Sudakov, and Verstra\"{e}te~\cite{MR2286514} -- they initiated a systematic study of the rainbow Tur\'{a}n number, where for a fixed $H$ and an integer $n$, the \emph{rainbow Tur\'{a}n number for $H$} is the maximum number of edges in a properly edge-colored graph on $n$ vertices which does not contain a rainbow $H$;
the other class is the edge-colored multi-graphs, which is equivalent to the graph system language in the abstract, and is the main object of study in this paper.
Although the two problems are different, a common scheme can be formulated as below.
A $k$-uniform hypergraph, \emph{k}-\emph{graph} for short, is a pair $H=(V, E)$, where $V$ is a finite set of vertices and $E\subseteq \binom{V}{k}$.
We identify a hypergraph $H$ with its edge set, writing $e\in H$ for $e\in E(H)$.
We write \emph{subgraph} instead of sub-$k$-graph or subhypergraph for brevity.

\begin{definition}
A $k$-graph system $\textbf{G}=\{G_1,\ldots,G_m\}$ is a collection of not necessarily distinct $k$-graphs on the same vertex set $V$.
Then a $k$-graph $H$ on $V$ is rainbow if there exists an injection $\varphi: E(H)\rightarrow[m]$ such that $e\in E(G_{\varphi(e)})$ for each $e\in E(H)$.
\end{definition}

Since $\varphi$ is an injection, it follows that all edges of $H$ are from different members of $\emph{\textbf{G}}$.
When $m=e(H)$, $\varphi$ is a bijection and thus $H$ contains exactly one edge from each $H_i$.
Note that each $G_i$ can be seen as the collection of edges with color $i$.
Given a $k$-graph system $\emph{\textbf{G}}$, the \emph{color set} of a graph $H$, denoted by $C(H)$, is the index set of all edges, i.e. $\{i:E(G_i)\cap E(H) \neq \emptyset\}$.
Note that if $|E(H)|=m$, then a rainbow $H$ consists of exactly one edge from each $G_i$.

As for the edge-colored multi-graphs, a recent breakthrough of Aharoni, DeVos, de la Maza, Montejano and \v{S}\'{a}mal~\cite{MR4125343} establishes a rainbow version of the Mantel's~\cite{mantel} theorem:  for $\textbf{\emph{G}}=\{G_1,G_2,G_3\}$ on the same $n$-vertex set, if $e(G_i)>\tau n^2$ for $i\in [3]$ where $\tau\approx 0.2557$, then $\textbf{\emph{G}}$ contains a rainbow triangle. Moreover, the constant $\tau$ is best possible.
Towards a better understanding of the rainbow structures, Aharoni~\cite{MR4125343} conjectured a rainbow version of the Dirac's~\cite{Dirac} theorem: for $|V|=n\geq 3$ and $\emph{\textbf{G}}=\{G_1,\dots, G_n\}$ on $V$, if $\delta(G_i)\geq n/2$ for each $i\in[n]$, then $\emph{\textbf{G}}$ contains a rainbow Hamilton cycle.
This was recently verified asymptotically by Cheng, Wang and Zhao~\cite{2019Rainbow}, and completely by Joos and Kim~\cite{MR4171383}.

A natural next step is to study rainbow analogues of graph factors, which we define now.
Given graphs $F$ and $G$, an $F$-\emph{tiling} is a set of vertex-disjoint copies of $F$ in $G$.
A \emph{perfect} $F$-$tiling$ (or an $F$-\emph{factor}) of $G$ is an $F$-tiling covering all the vertices of $G$.
Finding sufficient conditions for the existence of an $F$-factor is one of the central areas of research in extremal graph theory.
The celebrated Hajnal--Szemer\'{e}di theorem reads as follows.

\begin{theorem}[Hajnal--Szemer\'{e}di~\cite{MR0297607}, Corr\'{a}di--Hajnal~\cite{MR200185} for $t=3$]
\label{2}
Every $n$-vertex graph $G$ with $n\in t\mathbb{N}$ and $\delta(G)\geq(1-\frac{1}{t})n$ has a  $K_t$-factor.
Moreover, the minimum degree condition is sharp.
\end{theorem}

A short and elegant proof was later given by Kierstead and Kostochka \cite{2015A}.
The minimum degree threshold forcing an $F$-factor for arbitrary $F$ was obtained by K\"uhn and Osthus~\cite{Daniela2006Critical,Daniela2009The}, improving earlier results of Alon and Yuster~\cite{1992AlmostH} and Koml\'os, S\'ark\"ozy and Szemer\'edi~\cite{Komlos2001Proof}.




\subsection{Our results}
In this paper, we study rainbow clique-factors under a few different contexts.
Our first result is an asymptotical version of the rainbow Hajnal--Szemer\'edi theorem.

\begin{theorem}\label{3}
For every $\varepsilon>0$ and $t\in \mathbb N$, there exists $n_0\in \mathbb{N}$ such that the following holds for all integers $n\geq n_0$ and $n\in t\mathbb{N}$.
Let $m=\frac{n}{t}\binom{t}{2}$ and $\textbf{G}=\{G_1,\ldots, G_m\}$ be an $n$-vertex graph system.
If $\delta(G_i)\geq(1-\frac{1}{t}+\varepsilon)n$ for each $i$, then $\textbf{G}$ contains a rainbow $K_t$-factor.
\end{theorem}

The minimum degree conditions are asymptotically best possible, as seen by setting all $G_i$ to be identical and then referring to the optimality of Theorem~\ref{2}.
In fact, this naive construction serves as a simple lower bound for all ``rainbow'' problems, certifying that the rainbow version is at least ``as hard as'' the single host graph version (although the aforementioned rainbow Mantel's theorem says that the rainbow version can be strictly ``harder'').

It is natural to seek analogues of the Hajnal--Szemer\'{e}di theorem in the digraph and oriented graph settings, where we consider factors of directed cliques, namely, tournaments.
We consider digraphs with no loops and at most one edge in each direction between every pair of vertices.
Let $T_k$ be the transitive tournament on $k$ vertices, where a \emph{transitive tournament} is an orientation of a complete graph $D$ with the property that if $xy$ and $yz$ are arcs in $D$ with $x\neq z$, then the arc $xz$ is also in $D$.
The minimum semi-degree $\delta^0(G)$ of a digraph $G$ is the minimum of its minimum out-degree
$\delta^+(G)$ and its minimum in-degree $\delta^-(G)$.

Czygrinow, DeBiasio, Kierstead and Molla \cite{MR3371502} proved that every digraph $G$ on $n$ vertices with $\delta^+(G)\geq(1-1/k)n$ contains a perfect $T_k$-tiling for integers $n,k$ with $k\mid n$.
Let $\mathcal{T}_k$ be the family of all tournaments on $k$ vertices, Treglown \cite{MR3406450} proved that given an integer $k\geq3$, there exists an $n_0 \in \mathbb{N}$ such that the following holds.
Suppose $T\in \mathcal{T}_k$, and $D$ is a digraph on $n\geq n_0$ vertices where $k\mid n$. If $\delta^0(D)\geq(1-1/k)n$, then there exists a $T$-factor.
For more results, see~\cite{MR3240133,MR3406450}.
In this paper, we prove the following extensions of Theorem~\ref{3} for digraphs.

\begin{theorem}\label{direct}
For every integer $k\ge 3$ and real $\varepsilon>0$, there exists $n_0\in \mathbb{N}$ such that the following holds for all integers $n\geq n_0$ and $n\in k\mathbb{N}$.
If $\textbf{D}=\{D_1,\ldots,D_m\}$, $m=\frac{n}{k}\binom{k}{2}$, is a collection of $n$-vertex digraphs on the same vertex set such that $\delta^+(D_i)\geq(1-\frac{1}{k}+\varepsilon)n$, then $\textbf{D}$ contains a rainbow $T_k$-factor.
\end{theorem}

\begin{theorem}\label{semi}
For every integer $k\ge 3$, $T\in \mathcal{T}_k$ and real $\varepsilon>0$, there exists $n_0 \in \mathbb{N}$ such that the following holds for all integers $n\geq n_0$ and $n\in k\mathbb{N}$.
If $\textbf{D}=\{D_1,\ldots,D_m\}$, $m=\frac{n}{k}\binom{k}{2}$, is a  collection of $n$-vertex digraphs on the same vertex set such that $\delta^0(D_i)\geq(1-\frac{1}{k}+\varepsilon)n$, then $\textbf{D}$ contains a rainbow $T$-factor.
\end{theorem}

We next discuss the partite setting.
Suppose $V_1,\ldots, V_k$ are disjoint vertex sets each of order $n$, and $G$ is a $k$-partite graph on vertex classes $V_1,\ldots, V_k$ (that is, $G$ is a graph on the vertex set $V_1\cup\cdots\cup V_k$ such that no edge of $G$ has both end vertices in the same class).
We define the partite minimum degree of $G$, denoted by $\delta'(G)$, to be the largest $m$ such that every vertex has at least $m$ neighbours in each part other than its own, i.e.
\[
\delta'(G):= \min \limits_{i\in[k]}\min \limits_{v\in V_i}\min \limits_{j\in[k]\setminus\{i\}}|N(v)\cap V_j|,
\]
where $N(v)$ denotes the neighbourhood of $v$.
Fischer \cite{MR1698745} conjectured that if $\delta'(G)\geq(1-1/k)n$, then $G$ has a $K_k$-factor.
Recently, an approximate version of this conjecture assuming the degree condition $\delta'(G)\geq(1-1/k+o(1))n$ was proved independently by Keevash and Mycroft \cite{MR3290271}, and by Lo and Markstr\"{o}m \cite{MR3002575} (a corrected exact version was given by Keevash and Mycroft~\cite{KEEVASH2015187}, in fact, they obtain a more general result for $r$-partite graphs with $r\ge k$).
We extend this approximate version to the rainbow setting.

\begin{theorem}\label{rpartite}
For every $\varepsilon>0$ and integer $k$, there exists $n_0\in \mathbb{N}$ such that the following holds for all integers $n\geq n_0$.
If $\textbf{G}=\{G_1,\ldots,G_{n\binom{k}{2}}\}$ is a collection of $k$-partite graphs with a common partition $V_1,\ldots,V_k$ each of size $n$ such that $\delta'(G_i)\geq(1-\frac{1}{k}+\varepsilon)n$, then $\textbf{G}$ contains a rainbow $K_k$-factor.
\end{theorem}

A \emph{matching} in $H$  is a collection of vertex-disjoint edges of $H$. A \emph{perfect matching} in $H$ is a matching that covers all vertices of $H$.
Given $d\in [k-1]$, the minimum $d$-degree of a $k$-graph $H$, denoted by $\delta_d(H)$, is defined as the minimum of $d_H(S)$ over all $d$-sets $S$ of $V(H)$ where $d_H(S)$ denotes the number of edges containing $S$.
Another main result of this paper is to settle the rainbow version of minimum $d$-degree-type results for perfect matchings in $k$-graphs for all $d\in [k-1]$, in a sense that the minimum $d$-degree condition which \emph{forces a perfect matching in a single $k$-graph is essentially sufficient to force a rainbow perfect matching in a $k$-graph system}.
Note that Joos and Kim \cite{MR4171383} proved that $\delta(G_i)\geq n/2$ guarantees a rainbow perfect matching in an $n$-vertex graph system.


It is well-known that perfect matchings are closely related to its fractional counterpart.
Given a $k$-graph $H$, a \emph{fractional matching} is a function $f: E(H)\to [0,1]$, subject to the requirement that $\sum_{e:v\in e}f(e)\le 1$, for every $v\in V(H)$. Furthermore, if equality holds for every $v\in V(H)$, then we call the fractional matching \emph{perfect}.
Denote the maximum size of a fractional matching of $H$ by $\nu^*(H)=\max _f \Sigma_{e\in E(H)} f(e)$.
Let $c_{k,d}$ be the minimum $d$-degree threshold for perfect fractional matchings in $k$-graphs, namely, for every $\varepsilon>0$ and sufficiently large $n\in \mathbb N$, every $n$-vertex $k$-graph $H$ with $\delta_d(H)\ge (c_{k,d}+\eps)\binom{n-d}{k-d}$ contains a perfect fractional matching.
It is known that~\cite{MR2915641} every $n$-vertex $k$-graph $H$ with $\delta_d(H)\ge (\max\{c_{k,d}, 1/2\}+o(1))\binom{n-d}{k-d}$ has a perfect matching, and this condition is asymptotically best possible.
However, determining the parameter $c_{k,d}$ is a major open problem in this field and we refer to~\cite{2018The} for related results and discussions.

\begin{theorem}
\label{main}
For every $\varepsilon>0$ and integer $d\in[k-1]$, there exists $n_0\in \mathbb N$, such that the following holds for all integers  $n\geq n_0$ and $n\in k\mathbb{N}$.
Every $n$-vertex $k$-graph system $\textbf{G}=\{G_1, \dots, G_{\frac{n}{k}}\}$ with $\delta_d(G_i)\geq (\max\{c_{k,d}, \frac{1}{2}\}+\varepsilon)\binom{n-d}{k-d}$ for each $i$ contains a rainbow perfect matching.
\end{theorem}

\noindent\textbf{Related work}
Aharoni and Howard \cite{Preprint} conjectured that given an $n$-vertex $k$-graph system $\textbf{\emph{G}}=\{G_1,\ldots,G_m\}$, if $e(G_i)>\max\{\binom{n}{k}-\binom{n-m+1}{k},\binom{km-1}{k}\}$ for $i\in[m]$, then $\textbf{\emph{G}}$ contains a rainbow matching of size $m$.
The conjecture is known for $n>3k^2m$ by a result of Huang, Loh and Sudakov \cite{2012The}, and for $m<n/(2k)$ and sufficiently large $n$ by a recent result of Lu, Wang and Yu \cite{2020A}.
For the case $k=3$, Lu, Yu and Yuan~\cite{2020B} showed that for sufficiently large $n\in3\mathbb{N}$, given a $3$-graph system $\textbf{\emph{G}}=\{G_1,\ldots,G_{n/3}\}$, if $\delta_1(G_i)>\binom{n-1}{2}-\binom{2n/3}{2}$ for $i\in[n/3]$, then $\textbf{\emph{G}}$ contains a rainbow perfect matching (note that the single host 3-graph case was proved by K\"{u}hn, Osthus, and Treglown \cite{K2013Matchings} and independently by Khan \cite{IMDADULLAH2013PERFECT}).

On a slightly different setup, Huang, Li and Wang \cite{MR3904820} obtained a generalization of the Erd\H{o}s Matching Conjecture to properly-colored $k$-graph systems $\emph{\textbf G}$ and verified it for $n\geq3k^2m$ where $\textbf{\emph{G}}=\{G_1,\ldots,G_m\}$ and each $G_i$ is an $n$-vertex properly colored $k$-graph.
For general $F$-factors, Coulson, Keevash, Perarnau and Yepremyan~\cite{MR4055023} proved that essentially the minimum $d$-degree threshold guaranteeing an $F$-factor in a single $k$-graph also forces a rainbow $F$-factor in any edge-coloring of $G$ that satisfies certain natural local conditions.
We refer the reader to \cite{rainbow,Howard2017A,MR3604112,MR1652837,New,LU,Hongliang2018ON,MR3818098} for more results.



\section{Proof ideas and a general framework for rainbow $F$-factors}
Our proof is under the framework of the absorption method, pioneered by R\"odl, Ruci\'nski and Szemer\'edi \cite{MR2195584}, which reduces the problem of finding a spanning subgraph to building an absorption structure and an almost spanning structure.
Tailored to our problem, the naive idea is to build a rainbow absorption structure and a rainbow almost $F$-factor.
Moreover, the rainbow absorption structure must be able to deal with (i.e., absorb) an arbitrary leftover of vertices, as well as \emph{a leftover of colors}.

\subsection{A general framework for rainbow $F$-factors}
To state our general theorem for rainbow $F$-factors, we need some general notation that captures all of our contexts.
We shall consider a \emph{directed $k$-graph} (D$k$-graph) $H$, with edge set $E(H)\subseteq \binom{V(H)}k\times \{+, -\}$, that is, each edge consists of $k$ vertices and a direction taken from $\{+, -\}$.
This way a directed ($2$-)graph can be recognized as a graph with an \emph{ordered} vertex set, and edges following (or against) the order of the enumeration are oriented by $+$ (or $-$).


Given a D$k$-graph $H=(V,E)$, for $E'\subseteq E$, we write $H[E']$ for the subgraph of $H$ with edge set $E'$ and vertex set $\cup_{e\in E'} E'$.
For $V'\subseteq V$, if $H'\subseteq H$ contains all edges of $H$ with vertices in $V'$, then $H'$ is an \emph{induced subgraph} of $H$.
Denote it by $H[V']$.
If there is an $F$-tiling in $H$ whose vertex set is $V'$, then we say that $V'$ spans an $F$-tiling.
Given another D$k$-graph $H_1=(V_1,E_1)$, we set $H\cup H_1:= (V\cup V_1, E\cup E_1)$.

Given a D$k$-graph $F$ with $b$ vertices and $f$ edges, a D$k$-graph system $\emph{\textbf{G}}=\{G_1,\ldots,G_{\frac{n}{b}f}\}$ on vertex $V$ and a subset $V'\subseteq V$. Let $\emph{\textbf{G}}[V']=\{G_1[V'], \ldots, G_{\frac{n}{b}f}[V']\}$ be the \emph{induced} D$k$-graph system on $V'$. If $|V'|\in b\mathbb{N}$ and there exists a rainbow perfect $F$-tiling inside $\emph{\textbf{G}}[V']$ whose color set is $C\subseteq [nf/b]$, then we say that $V'$ \emph{spans} a rainbow $F$-tiling in $\emph{\textbf{G}}$ with color set $C$.
Let $A_1$ and $A_2$ be two rainbow D$k$-graphs in $\emph{\textbf{G}}$ with color set $C_1$ and $C_2$ respectively, we set $A_1\cup A_2$ to be a D$k$-graph with vertex set $V(A_1)\cup V(A_2)$, edge set $E(A_1)\cup E(A_2)$ and color set $C_1\cup C_2$.
The following general minimum degree condition captures all of our contexts.

\begin{definition}
Let $H$ be a D$k$-graph and $d\in [k-1]$, a minimum $d$-degree $\delta^*_d(H)$ corresponding to a certain (implicit) degree rule can be defined as follows.
There exists $\ell\in \mathbb{N}$ such that for any $d$-set $S\subseteq V(H)$, let
\[
N^*(S):=\{E_1,\ldots,E_{\ell}\},
\]
where $E_i$ is a set of edges of $H$ that each contains $S$ as a subset.
Let $\deg^*(S):=\min_{i\in [\ell]}|E_i|$ and $\delta^*_d(H):=\min_{S\subseteq\binom{V}{d}}\deg^*(S)$.
\end{definition}

We list all the instances of minimum degrees used in this paper.
\begin{itemize}
  \item When $H$ is a $k$-graph, we can take $\ell=1$ and $E_1$ as all edges containing $S$, thus $\deg^*(S)=|E_1|$ and $\delta^*_d(H)$ represents the standard minimum $d$-degree of $H$.
  \item When $H$ is a $k$-partite 2-graph, we can take $\ell=k-1$ where $E_i$ consists of the edges from $S=\{v\}$ to $V_i$ for $i\in[k]\setminus\{j:v\in V_j\}$, thus $\deg^*(S)=\min\{|E_i|:i\in[k]\setminus\{j:v\in V_j\}\}$ and $\delta^*_1(H)$ represents the minimum partite degree of $H$.
  \item When $H$ is a directed graph and $S=\{v\}$ for a given vertex $v$, we can take $\ell=1$ and take $E_1$ as the sets of out(in)-edges, thus,  $\deg^*(S)=|E_1|$ and $\delta^*_1(H)$ represents the minimum out(in)-degree of $H$.
  \item When $H$ is a directed graph and $S=\{v\}$ for a given vertex $v$, we can take $\ell=2$ and take $E_1, E_2$ as the sets of in-edges and out-edges respectively, thus, $\deg^*(S)=\min\{|E_1|,|E_2|\}$ and $\delta^*_1(H)$ represents the minimum semi-degree of $H$.
\end{itemize}

Throughout the rest of this paper, let $F$ be a D$k$-graph with $b$ vertices and $f$ edges.
We first define an absorber without colors.
Given a set $B$ of $b$ vertices, a D$k$-graph $A^0=A_1^0\cup A_2^0$ is called an \emph{$F$-absorber for $B$} if
\begin{itemize}
  \item $V(A^0)=B\dot\cup L$\footnote{As usual, $A\dot\cup B$ denotes the disjoint union of $A$ and $B$.},
  \item $A_1^0$ is an $F$-factor on $L$ and $A_2^0$ is an $F$-factor on $B\cup L$.
\end{itemize}
Note that $|V(A^0)|$ is always a constant in this paper.
Naturally, we give the definition of rainbow $F$-absorber as follows.

\begin{definition}[Rainbow $F$-absorber]
Let $\textbf G=\{G_1,\dots, G_{nf/b}\}$
be a D$k$-graph system on $V$ and $F$ be a D$k$-graph with $b$ vertices and $f$ edges.
For every $b$-set $B$ in $V$ and every $f$-set $C$ in $[nf/b]$, $A=A_1\cup A_2$ is called a rainbow $F$-absorber for $(B,C)$ if
\begin{itemize}
  \item $V(A)=B\dot\cup L$,
  \item $A_1$ is a rainbow $F$-factor on $L$ with color set $C_1$ and $A_2$ is a rainbow $F$-factor on $B\cup L$ with color set $C_1\cup C$.
\end{itemize}
\end{definition}

A rainbow $F$-absorber for $(B, C)$ works in the following way:
$A_1$ is a rainbow $F$-tiling with color set $C_1$,
thus, if the vertices of $B$ and the colors $C$ are available, then we can switch to $A_2$ and get a larger rainbow $F$-tiling.
For example, Let $\emph{\textbf{G}}=\{G_1,\ldots,G_{n}\}$ be a graph system on $n$-vertex set $V$.
For any $B=\{v_1,v_2,v_3\}$ and $C=\{i,j,k\}$, we construct a rainbow triangle-absorber $A$ for $(B,C)$ where $V(A)=B\cup \{v_4,\ldots,v_{12}\}$.
$\{v_4v_7v_8,v_5v_9v_{10},v_6v_{11}v_{12}\}$ is a family of rainbow triangles with color set $C_1$, which serves as $A_1$.
$\{v_1v_7v_8,v_2v_9v_{10},v_3v_{11}v_{12},v_4v_5v_6\}$ is a family of rainbow triangles with color set $C_1\cup\{i,j,k\}$,  which serves as $A_2$.
\begin{figure}[htb]
\centering
\includegraphics[width=12cm]{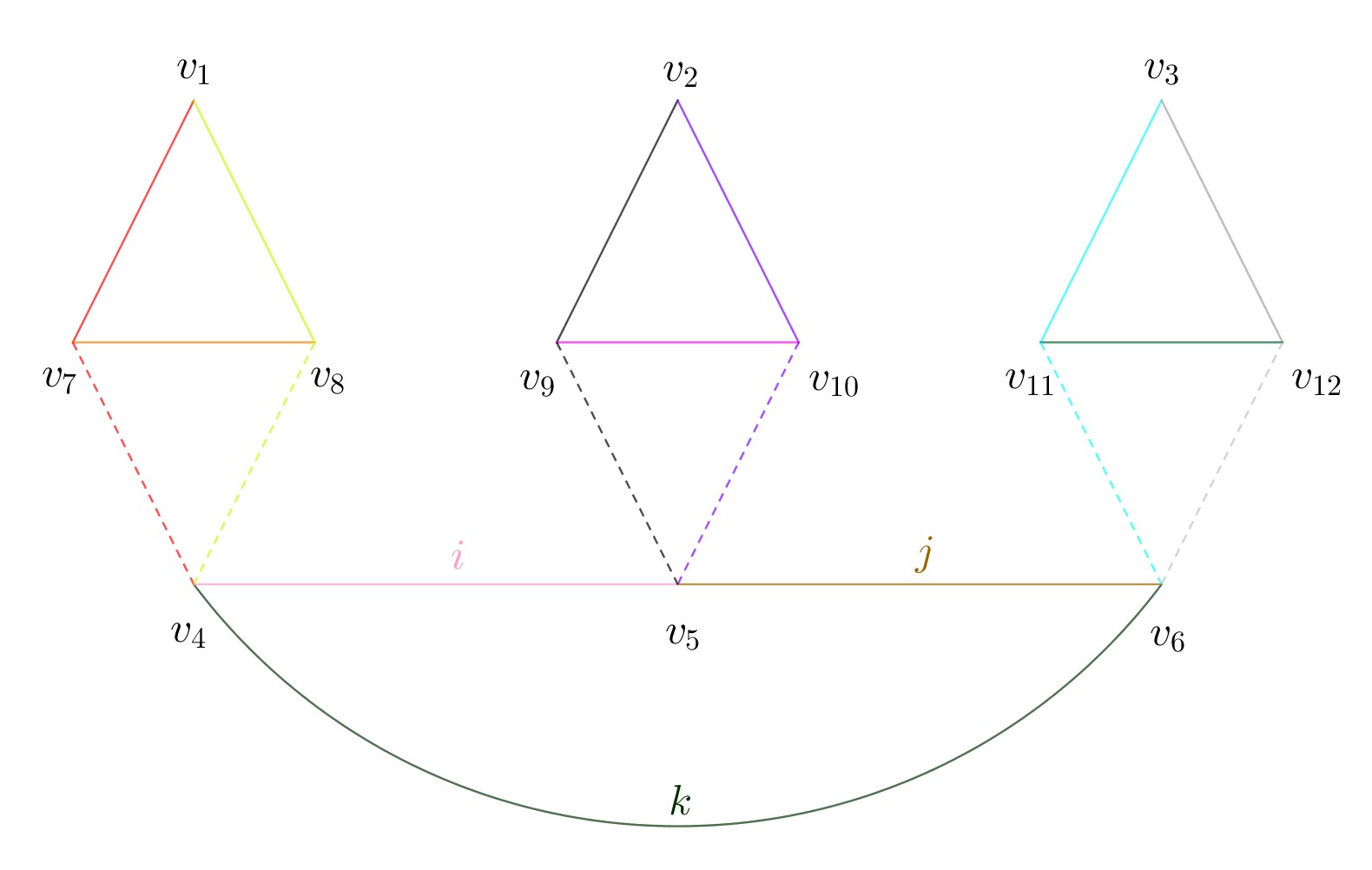}\\
\caption{A rainbow triangle-absorber $A$ for $(B,C)$}
\label{absorber3}
\end{figure}

Now we introduce one of the main parameters $c_{d,F}^{\rm abs,*}$.
Roughly speaking, it is the minimum degree threshold such that all $b$-sets are contained in many rainbow $F$-absorbers.

\begin{definition}[$c_{d,F}^{\rm abs,*}$: Rainbow absorption threshold]\label{abs}
Fix an $F$-absorber $A^0=A^0_1\cup A^0_2$ and let $m$ be the number of vertex disjoint copies of $F$ in $A^0_2$.
Let $c_{d,F,A^0}\in (0,1)$ be the infimum of reals $c>0$ such that for every $\eps>0$ there exists $\eps'>0$ such that the following holds for sufficiently large $n\in \mathbb{N}$ where $d\in[k-1]$.
Let $\textbf G=\{G_1,\dots, G_{nf/b}\}$ be an $n$-vertex D$k$-graph system on $V$.
If $\delta_d^*(G_i)\ge (c+\eps)\binom{n-d}{k-d}$ for $i\in[nf/b]$, then for every $b$-set $B$ in $V$ and every $f$-set $C$ in $[nf/b]$ with the form $[(i-1)f,if]$ for some $i\in[n/b]$, there are at least $\eps' n^{(m-1)(b+1)}$ rainbow $F$-absorbers $A$ with color set $C(A_1^0)\cup C$ whose underlying graph is isomorphic to $A^0$ such that $C(A_1^0)=[(i_1-1)f+1,i_1f]\cup[(i_2-1)f+1,i_2f]\cup\cdots\cup[(i_{m-1}-1)f+1,i_{m-1}f]$ where $i_j\in[n/b]$ for each $j\in[m-1]$ and $i_{j_1}\neq i_{j_2}$ for distinct $j_1,j_2\in[m-1]$.
Let $c_{d,F}^{\rm abs,*}:=\inf c_{d,F,A^0}$ where the infimum is over all $F$-absorbers $A^0$.
\end{definition}


We next define a threshold parameter for the rainbow almost $F$-factor in a similar fashion.
We use the following auxiliary $b$-graph $H_F$.
Given a D$k$-graph $F$ with $b$ vertices and $f$ edges, and an $n$-vertex D$k$-graph system $\emph{\textbf H}=\{H_1,\dots, H_{f}\}$ on $V$, let $H_F$ be the (undirected) $b$-graph with vertex set $V(H_F)=V$ and edge set $E(H_F)=\{V(F'):F'\ {\rm is\ a\ rainbow\ copy\ of}\ F\ {\rm with\ color\ set}\ [f]\}$.

\begin{definition}[$c_{d,F}^{\rm cov,*}$: Rainbow almost $F$-factor threshold]\label{cov}
Let $c_{d,F}^{\rm cov,*}\in (0,1)$ be the infimum of reals $c>0$ such that  for every $\eps>0$, the following holds for sufficiently large $n\in \mathbb{N}$.
Let $\textbf H=\{H_1,\dots, H_{f}\}$ be an $n$-vertex D$k$-graph system.
If $\delta_d^*(H_i)\ge (c+\eps)\binom{n-d}{k-d}$ for every $i\in[f]$, then the $b$-graph $H_F$ has a perfect fractional matching.
\end{definition}

The property of ``having a perfect fractional matching'' is required in $H_F$, which is a single host graph.
This definition (and our proofs supporting it) establishes a close relation between the rainbow $F$-factor problem and the classical $F$-factor problem with no colors.

Now we are ready to state our general result on rainbow $F$-factors.
\begin{theorem}
\label{general}
Let $F$ be a D$k$-graph with $b$ vertices and $f$ edges.
For any $\eps>0$ and integer $d\in [k-1]$, the following holds for sufficiently large $n\in b\mathbb{N}$.
Let $\textbf G=\{G_1,\dots, G_{nf/b}\}$ be an $n$-vertex D$k$-graph system on $V$.
If $\delta_d^*(G_i)\ge (\max\{c_{d, F}^{\rm abs,*}, c_{d, F}^{\rm cov,*}\}+\eps)\binom{n-d}{k-d}$ for $i\in[nf/b]$, then $\textbf G$ contains a rainbow $F$-factor.
\end{theorem}

%

Theorem~\ref{general} reduces the rainbow $F$-factor problem to two subproblems, namely, the enumeration of rainbow $F$-absorbers and the study of perfect fractional matchings in $H_F$.
In our proofs of Theorems~\ref{3}--\ref{main}, the first subproblem is done by greedy constructions of the $K_t$-absorbers with the minimum degree condition.
Note that the second subproblem is trivial by definition for Theorem~\ref{main}.
For Theorems~\ref{3}--\ref{rpartite}, we achieve it by converting the problem to the setting of complexes (downward-closed hypergraphs) and then applying a result of Keevash and Mycroft \cite{MR3290271} on perfect fractional matchings, which is a nice application of the Farkas' Lemma for linear programming.

To conclude, we remark that the main benefit from Theorem~\ref{general} is that both of these two subproblems only concern \emph{finitely many} colors.
From this aspect, Theorem~\ref{general} irons out significant difficulties on the rainbow spanning structure problem due to an unbound number of colors.
Thus, it is likely that Theorem~\ref{general} will find more applications in this area.


\section{Notation and Preliminary}

For a hypergraph $H$, the 2-degree of a pair of vertices is the number of edges containing this pair and $\Delta_2(H)$ denotes the maximum 2-degree in $H$.
For reals $a, b$ and $c$, we write $a = (1\pm b)c$ for
$(1-b)c \leq a\leq(1+b)c$.
We need the following result which was attributed to Pippenger \cite{pippenger}(see Theorem 4.7.1 in \cite{Probabilistic}), following Frankl and R\"{o}dl~\cite{frankl1985near}. A \emph{cover} in a hypergraph $H$ is a set of edges such that each vertex of $H$ is in at least one edge of the set.

\begin{lemma}[\cite{pippenger}]\label{pipp}
For every integer $k\geq2$, $r\geq1$ and $a>0$, there exist $\gamma=\gamma(k,r,a)>0$ and $d_0=d_0(d,r,a)$ such that the following holds for every $n\in \mathbb{N}$ and $D\geq d_0$.
Every $k$-graph $H=(V,E)$ on $V$ of $n$ vertices in which all vertices have positive degrees and which satisfies the following conditions:
\begin{itemize}
  \item For all vertices $x\in V$ but at most $\gamma n$ of them, $d_H(x)=(1\pm\gamma)D$.
  \item For all $x\in V$, $d_H(x)<rD$.
  \item $\Delta_2(H)<\gamma D$.
\end{itemize}
contains a cover of at most $(1+a)(n/k)$ edges.
\end{lemma}




The following well-known concentration results, i.e. Chernoff bounds, can be found in Appendix A in \cite{Probabilistic} and Theorem 2.8, inequality (2.9) and (2.11) in \cite{Random}.
Denote a binomial random variable with parameters $n$ and $p$ by $Bi(n,p)$.
Bernoulli distribution is the discrete probability distribution of
a random variable which takes the value 1 with probability $p$ and the value 0 with probability $1-p$.

\begin{lemma}[Chernoff Inequality for small deviation]\label{chernoff1}
If $X=\sum_{i=1}^nX_i$ where $X_1,\ldots,X_n$ are mutually independent random variables, each $X_i$ has Bernoulli distribution with expectation $p_i$ and $\alpha\leq3/2$, then
\[
\mathbb{P}[|X-\mathbb{E}[X]|\geq\alpha\mathbb{E}[X]]\leq2e^{-\frac{\alpha^2}{3}\mathbb{E}[X]}.
\]
In particular, when $X\sim Bi(n,p)$ and $\lambda<\frac{3}{2}np$, then
\[
\mathbb{P}[|X-np|\geq\lambda]\leq e^{-\Omega(\lambda^2/(np))}.
\]
\end{lemma}
\begin{lemma}[Chernoff Inequality for large deviation]\label{Chernoff2}
If $X=\sum_{i=1}^nX_i$ where $X_1,\ldots,X_n$ are mutually independent random variables, each random variable $X_i$ has Bernoulli distribution with expectation $p_i$ and $x\geq7\mathbb{E}[X]$, then
\[
\mathbb{P}[X\geq x]\leq e^{-x}.
\]
\end{lemma}

We also need the Janson's inequality~\cite{Probabilistic} to provide an exponential upper bound for the lower tail of a sum of dependent zero-one random variables.
\begin{lemma}[Theorem 8.7.2 in \cite{Probabilistic}]\label{Janson}
Let $\Gamma$ be a finite set and $p_i\in[0,1]$ be a real for $i\in\Gamma$. Let $\Gamma_p$ be a random subset of $\Gamma$ such that the elements are chosen independently with $\mathbb{P}[i\in\Gamma_p]=p_i$ for $i\in\Gamma$. Let $M$ be a family of subsets of $\Gamma$. For every $A_i\in M$, let $I_{A_i}=1$ if $A_i\subseteq\Gamma_p$ and 0 otherwise.
Let $B_i$ be the event that $A_i\subseteq\Gamma_p$.
For $A_i, A_j\in M$, we write $i\sim j$ if $B_i$ and
$B_j$ are not pairwise independent, in other words, $A_i\cap A_j\neq\emptyset$.
Define $X=\Sigma_{A_i\in M}I_{A_i}$, $\lambda=\mathbb{E}[X]$, $\Delta=\sum\limits_{i\sim j}\mathbb{P}[B_i\wedge B_j]$, then
\[
\mathbb{P}[X\leq(1-\gamma)\lambda]<e^{-\gamma^2\lambda/[2+(\Delta/\lambda)]}.
\]
\end{lemma}

\section{Rainbow Absorption Method}
\subsection{Rainbow Absorption Lemma}
In this section we prove a rainbow version of the absorption
lemma via probabilistic method.
The only difference is that we use the following auxiliary hypergraph which makes it applicable to the rainbow setting.

\begin{definition}\label{fbg}
We call a hypergraph $H$ a $(1,b)$-graph, if $V(H)$ can be partitioned into $A\cup B$ and $E(H)$ is a family of $(1+b)$-sets each of which contains exactly one vertex in $A$ and $b$ vertices in $B$.
\end{definition}


For a $(1,b)$-graph $H$ with partition $A\dot\cup B$, a $(1,d)$-subset $D$ of $V(H)$ is a $(d+1)$-tuple where $|D\cap A|=1$ and $|D\cap B|=d$.
A $(1,b)$-graph $H$ with partition classes $A,B$ is \emph{balanced }if $b|A|=|B|$.
We say that a set $S\subseteq V(H)$ is balanced if $b|S\cap A|=|S\cap B|$.

Given an $n$-vertex D$k$-graph system $\emph{\textbf{G}}=\{G_1,\ldots,G_{nf/b}\}$ on $V$,
we first construct a sequence of hypergraphs $H_{F_1}',\ldots,H_{F_{n/b}}'$, each of which is a $b$-graph with vertex set $V(H_{F_i}')=V$ and edge set $E(H_{F_i}')=\{e\in\binom{V}{b}:$ $e$ spans a rainbow copy of $F$ with color set $I_i=[(i-1)f+1,if]$$\}$.
We define an auxiliary $(1,b)$-graph $H_{\emph{\textbf{G}}}$ of $\emph{\textbf{G}}$ as follows.
\begin{definition}
Let $H_{\textbf{G}}$ be an auxiliary $(1,b)$-graph of $\textbf{G}$ with vertex set $V'=[n/b]\cup V$ and edge set $\{\{i\}\cup e: i\in[n/b], e\in H_{F_i}'\}$.
\end{definition}
For any edge $e\in E(H_{\emph{\textbf{G}}})$,
If $A\subseteq V(H_{\emph{\textbf{G}}})$ and $|A|$ is divisible by $b+1$, then $A\in\binom{(b+1)n}{a}$
is an \emph{absorber} for $e$ if $e\subseteq A$, there is a perfect matching in $H_{\emph{\textbf{G}}}[A]$ and there is a perfect matching in $H_{\emph{\textbf{G}}}[A\setminus e]$.
Let $\mathcal{L}(e)$ denote the set of absorbers for $e$ in $H_{\emph{\textbf{G}}}$.

\begin{lemma}[Rainbow Absorption Lemma]
\label{absorption}
Let $F$ be a D$k$-graph with $b$ vertices and $f$ edges and $A^0$ be a rainbow $F$-absorber.
The maximum vertex-disjoint copies of $F$ of $A^0$ is $m$.
For every $\varepsilon>0$, there exist $\gamma,\gamma_1$ and $n_0$ such that the following holds for all integers $n\geq n_0$.
Suppose that $\textbf{G}=\{G_1,\ldots,G_{\frac{n}{b}f}\}$ is an $n$-vertex D$k$-graph system on $V$ and $\delta_d^*(G_i)\geq(c_{d,F}^{\rm abs,*}+\varepsilon)\binom{n-d}{k-d}$ and $H_{\textbf{G}}$ is the auxiliary $(1,b)$-graph of $\textbf{G}$, then there exists a matching $M$ in $H_{\textbf{G}}$ with size at most $2\gamma (m-1)n$ such that for every balanced set $U\subseteq \left([n/b]\cup V\right)\setminus V(M)$ of size at most $\gamma_1n$,
$V(M)\cup U$ spans a matching in $H_{\textbf{G}}$.
\end{lemma}
\begin{proof}
Let $1/n\ll\gamma_1\ll\alpha\ll\gamma\ll\varepsilon'\ll\varepsilon$.
Note that a matching of size $m$ in $H_{\emph{\textbf{G}}}$ corresponds to a rainbow $F$-absorber in $\textbf{\emph{G}}$.
Choose a family $\mathcal{F}$ of matchings of size $m-1$ from $H_{\emph{\textbf{G}}}$ by including each matching of size $m-1$ independently at random with probability
\[
p=\gamma/n^{(m-1)(b+1)-1}.
\]
Note that $|\mathcal{F}|, |\mathcal{L}(e)\cap\mathcal{F}|$ are binomial random variables with expectations
\[
\mathbb{E}|\mathcal{F}|\leq \gamma n\ {\rm and}
\]
\[
\mathbb{E}|\mathcal{L}(e)\cap\mathcal{F}|\geq\gamma\varepsilon'n {\rm\ for\ any}\ e\in E(H_{\emph{\textbf{G}}}).
\]
The latter inequality holds since for any edge $e$ of $H_{\emph{\textbf{G}}}$, $|\mathcal{L}(e)|\geq\varepsilon'n^{(m-1)(b+1)}$ by the minimum degree assumption and Definition~\ref{abs}.
By Lemma~\ref{chernoff1}, with probability $1-o(1)$, the family $\mathcal{F}$ satisfies the following properties.
\begin{enumerate}
  \item $|\mathcal{F}|\leq2\mathbb{E}|\mathcal{F}|\leq 2\gamma n$,
  \item $|\mathcal{L}(e)\cap\mathcal{F}|\geq\frac{1}{2}\mathbb{E}|\mathcal{L}(e)\cap\mathcal{F}|\geq\frac{1}{2}\gamma\varepsilon'n$ for any $e\in E(H_{\emph{\textbf{G}}})$.
\end{enumerate}
Moreover, we can also bound the expected number of pairs of intersecting members of $\mathcal{F}$ by
\[
n^{(m-1)(b+1)}(m-1)^2(b+1)^2n^{(m-1)(b+1)-1}p^2\leq \frac{1}{8}\gamma\varepsilon'n.
\]
Thus, by Markov's inequality, we derive that with probability at least $1/2$, $\mathcal{F}$ contains at most $\frac{1}{4}\gamma\varepsilon'n$ pairs of intersecting members of $\mathcal{F}$.
Remove one member from each of the intersecting pairs in $\mathcal{F}$.
Thus, the resulting family, say $\mathcal{F}'$, consists of pairwise disjoint matchings of size $m-1$ that satisfies
\begin{enumerate}
  \item $|\mathcal{F}'|\leq 2\gamma n$,
  \item $|\mathcal{L}(e)\cap\mathcal{F}|\geq\frac{1}{2}\gamma\varepsilon'n-\frac{1}{4}\gamma\varepsilon'n\geq\alpha n$ for any $e\in E(H_{\emph{\textbf{G}}})$.
\end{enumerate}
Therefore, the union of members in $\mathcal{F}'$ is a matching in $H_{\emph{\textbf{G}}}$ of size at most $2\gamma(m-1)n$ and can (greedily) absorb a balanced set $U$ of size at most $\gamma_1n$ since $\gamma_1\ll\alpha$.
\end{proof}
\subsection{Enumeration of Rainbow Absorbers}
In this section, we give two examples of rainbow absorbers.
The first one is as follows.
Recall that $T_k$ is the transitive tournament on $k$ vertices.
\subsubsection{Rainbow $T_k$-absorber}
\label{Rainbow Absorption Method}
\begin{figure}[htb]
\centering
\includegraphics[width=12cm]{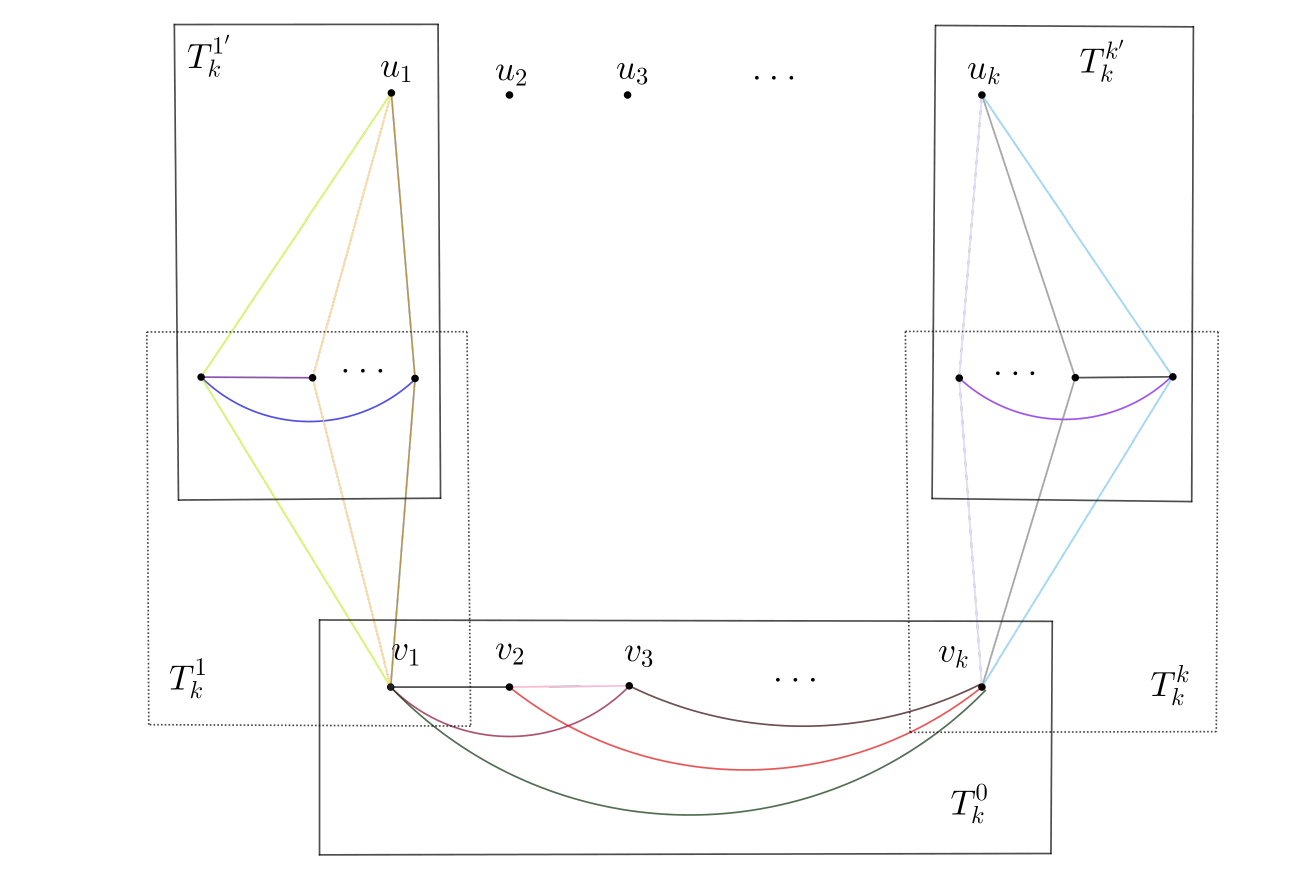}\\
\caption{Illustration of the rainbow absorbers (with directions omitted)}
\label{absorber}
\end{figure}

For the proof of Theorem~\ref{direct}, we show that
\begin{equation}
\label{absdi}
c_{1,T_k}^{\rm abs,+}\leq1-\frac{1}{k}.
\end{equation}
For any $k$-set $S=\{u_1,u_2,\ldots,u_k\}$ in $V$ and every $\binom{k}{2}$-set $C=[(j-1)\binom{k}{2}+1,j\binom{k}{2}]$ where $j\in[\frac{n}{k}]$, we define a rainbow $T_k$-absorber $A$ for $(S,C)$ as follows (see Figure \ref{absorber}).
\begin{itemize}
  \item $A_1=A-S=\{T_k^1,\ldots,T_k^k\}$ is a rainbow $T_k$-tiling with $C(T_k^i)=[(j_i-1)\binom{k}{2}+1,j_i\binom{k}{2}]$ where $j_i\in[n/k]$ for $i\in[k]$ and $A_2=\{T_k^{0},T_k^{1'},\ldots,T_k^{k'}\}$ is a rainbow $T_k$-tiling with $C(T_k^{0})=C$ and $C(T_k^{i'})=[(j_i-1)\binom{k}{2}+1,j_i\binom{k}{2}]$ where $j_i\in[n/k]$ for $i\in[k]$ .
  \item We can choose a rainbow $T_k^{0}$ that is isomorphic to $T_k$ with color set $C$ such that $V(T_k^{0})=\{v_1,v_2,\ldots,v_k\}$ where $v_i\in V(T_k^i)$ and $V(T_k^{i'})=V(T_k^i)\backslash \{v_i\}\cup\{u_i\}, i\in[k]$.
  \item $c(u_ix)=c(v_ix)$ for each $x\in V(A_1\backslash V(T_k^{0}))$, and $c(xy)$ in $T_k^i$ is the same as $c(xy)$ in $T_k^{i'}$ for $i\in[k]$ and $x,y\in V(A_1\backslash V(T_k^{0}))$.
\end{itemize}

Suppose $\delta(D_i^+)\ge(1-1/k+\varepsilon)n$ for $i\in[\frac{n}{k}\binom{k}{2}]$.
For any $k$-set $S$ in $V$ and every $\binom{k}{2}$-set $C=[(j-1)\binom{k}{2}+1,j\binom{k}{2}]$ where $j\in[\frac{n}{k}]$, we denote the family of rainbow $T_k$-absorbers for $(S,C)$ by $\mathcal{A}(S,C)$.
\begin{claim}\label{A}
For any $k$-set $S=\{u_1,u_2,\ldots,u_k\}$ in $V$ and every $\binom{k}{2}$-set $C=[(j-1)\binom{k}{2}+1,j\binom{k}{2}]$ where $j\in[\frac{n}{k}]$, we have $|\mathcal{A}(S,C)|\geq\varepsilon^{k^2+k}n^{k^2+k}$.
\end{claim}

\begin{proof}
Fixing a $k$-set $S=\{u_1,u_2,\ldots,u_k\}$ in $V$ and a $\binom{k}{2}$-set $C=[(j-1)\binom{k}{2}+1,j\binom{k}{2}]$ for some $j\in[\frac{n}{k}]$, we construct rainbow absorbers for $(S,C)$.
We choose $[(j_{i}-1)\binom{k}{2}+1,j_i\binom{k}{2}]$ for $i\in[k]$ arbitrarily, there are $(\frac{n}{k}-1)(\frac{n}{k}-2)\cdots(\frac{n}{k}-k)\geq\varepsilon^kn^k$ choices.
Next, we choose a rainbow $T_k^0$ with color set $C$. Due to the minimum out-degree of $D_i$, the number of choices for $T_k^0$ is at least
\[
(n-k)\left((1-\frac{1}{k}+\varepsilon)n-(k+1)\right) \cdots \left(\frac{1}{k}+(k-1)\varepsilon n-(2k-1)\right)\geq \varepsilon^kn^k.
\]

Now we fix one such $U=\{v_1,v_2,\ldots,v_k\}$.
For each $i\in[k]$ and each pair $\{u_i,v_i\}$, suppose we succeed in choosing a set $S_i$ such that $S_i$ is disjoint to $W_{i-1} = \cup_{j\in[i-1]}S_j\cup S\cup U$, then $V(T_k^{i'})=S_i\cup\{u_i\}$ spans a rainbow $T_k$ in $\emph{\textbf{D}}$ with color set $[(j_i-1)\binom{k}{2}+1,j_i\binom{k}{2}]$ while so does $V(T_k^i)=S_i\cup\{v_i\}$.

For the first vertex in $S_1$, the number of choices is at least $(1-\frac{2}{k}+2\varepsilon)n-2k$, for the last vertex in $S_1$,  the number of choices is in  $k\varepsilon n-(3k-1)$. Since $\frac{1}{n}\ll \varepsilon$, the number of choices for $S_1$ is at least
\[
\left((1-\frac{2}{k}+2\varepsilon)n-2k\right) \cdots (k\varepsilon n-(3k-1))\geq\varepsilon^{k-1}n^{k-1}.
\]

Similarly, the minimum out-degree implies that for $i\in[2,k]$ there are at least $\varepsilon^{k-1}n^{k-1}$ choices for $S_i$ and in total we obtain $\varepsilon^{k^2+k}n^{k^2+k}$ rainbow $T_k$-absorbers for $S$.
\end{proof}

\subsubsection{Rainbow edge-absorber}
\label{rainbow edge absorber}
For the proof of Theorem~\ref{main}, we show that
\begin{equation}
\label{absmatching}
c_{d,F}^{\rm{abs},d}\leq\frac{1}{2}.
\end{equation}
Given an $n$-vertex $k$-graph system $\emph{\textbf{G}}$ on $V$ with $\delta_d(G_i)\geq(\frac{1}{2}+\varepsilon)\binom{n-d}{k-d}$ for $i\in[n/k]$,
we first construct a $(1,k)$-graph $H_{\emph{\textbf{G}}}$ with vertex set $[n/k]\cup V$ and edge set $\{\{i\}\cup e:e\in H_i,i\in[n/k]\}$.
Next, we construct a specific rainbow edge-absorber.
For any $k$-set $T=\{v_1,\ldots,v_k\}$ in $V$ and every color $c_1\in[n/k]$, we give a rainbow absorber $A=A_1\cup A_2$ for $(T,c_1)$ as follows.
\begin{itemize}
  \item $A_1=\{M_2,\ldots,M_k\}$ is a set of $k-1$ disjoint edges in $H_{\emph{\textbf{G}}}$ where $c_i\in M_i(i\in[2,k])$.
  \item There is a vertex $u_i(i\in[2,k])$ from each $V(M_i)$ such that $\{u_2,\ldots,u_k, v_1, c_1\}\in E(H_{\emph{\textbf{G}}})$ and $(V(M_i)\setminus \{u_i\})\cup \{v_i\}\in E(H_{\emph{\textbf{G}}})$ for $i\in[2,k]$.
      Let $A_2$ be $\{\{u_2,\ldots,u_k, v_1, c_1\},(V(M_2)\setminus \{u_2\})\cup \{v_2\},\ldots,(V(M_k)\setminus \{u_k\})\cup \{v_k\}\}$.
\end{itemize}

For any $k$-set $T$ in $V$ and every color $c_1\in[n/k]$, we denote the family of such rainbow edge-absorbers for $(T,c_1)$ by $\mathcal{A}(T,c_1)$.

\begin{claim}\label{cl1}
$|\mathcal{A}(T,c_1)|\geq \varepsilon^{2k-2}n^{k-1} \binom{n-1}{k-1}^k/2$.
\end{claim}

\begin{proof}
Fix $c_1\in[n/k]$ and $T=\{v_1,\ldots,v_k\}\subseteq V$.
Choose $(c_2,\ldots,c_k)$ arbitrarily from $[n/k]$ and there are at least $(\frac{n}{k}-1)\cdots(\frac{n}{k}-(k-1))\geq\varepsilon^{k-1}n^{k-1}$ choices.
Fix such $(c_2,\ldots,c_k)$.
Next, we construct $M_2,\ldots,M_k$ and
note that there are at most $(k-1)\binom{n-1}{k-2}\leq \varepsilon\binom{n-1}{k-1}$ edges which contain $c_1, v_1$ and $v_j$ for some $j\in[2,k]$.
Due to the minimum degree assumption, there are at least $\frac{1}{2}\binom{n-1}{k-1}$ edges containing $v_1$ and $c_1$ but none of $v_2,\ldots,v_k$.
We fix such one edge $\{c_1,v_1,u_2,\ldots,u_k\}$ and set $U_1=\{u_2,\ldots,u_k\}$.
For each $i\in[2,k]$ and each pair $\{u_i,v_i\}$, suppose we succeed in choosing a set $U_i$ such that $U_i$ is disjoint with $W_{i-1}=\cup_{j\in[i-1]}U_j\cup T$ and both $U_i\cup\{u_i,c_i\}$ and $U_i\cup\{v_i, c_i\}$ are edges in $\tilde{H}$, then for a fixed $i\in[2,k]$, we call such a choice $U_i$ good.

Note that in each step $i\in[2,k]$, there are $k+(i-1)(k-1)\leq k^2$ vertices in $W_{i-1}$, thus the number of edges with color $c_i$ intersecting $u_i$ and at least one other vertex in $W_{i-1}$ is at most $k^2\binom{n-1}{k-2}$.
So the minimum degree assumption implies that for each $i\in[2,k]$, there are at least $2\varepsilon\binom{n-1}{k-1}-2k^2\binom{n-1}{k-2}\geq \varepsilon\binom{n-1}{k-1}$ choices for $U_i$ and in total we obtain $\varepsilon^{2k-2}n^{k-1} \binom{n-1}{k-1}^k/2$ rainbow absorbers for $(T,c_1)$.
\end{proof}

\section{Rainbow Almost Cover}

The goal of this section is to prove the following lemma, an important component of the proof of Theorem~\ref{general}.

\begin{lemma}[Rainbow Almost Cover Lemma]\label{rainbow cov}
Let $F$ be a D$k$-graph with $b$ vertices and $f$ edges.
For every $\varepsilon, \phi>0$ and integer $d\in[k-1]$, the following holds for sufficiently large $n\in b\mathbb{N}$.
Suppose that $\textbf{G}=\{G_1,\ldots,G_{nf/b}\}$ is an $n$-vertex D$k$-graph system on $V$ such that $\delta_d^*(G_i)\geq (c_{d,F}^{\rm cov,*}+\varepsilon)\binom{n-d}{k-d}$ for $i\in[nf/b]$, then $\textbf{G}$ contains a rainbow $F$-tiling covering all but at most $\phi n$ vertices.
 \end{lemma}

For a $k$-graph $H$, a \emph{fractional cover} is a function $\omega:V(H)\to [0,1]$, subject to the requirement $\sum_{v:v\in e}\omega(v)\geq1$ for every $e\in E(H)$.
Denote the minimum fractional cover size by $\tau^*(H)=\min_\omega\Sigma_{v\in V(H)}\omega(v)$.
The conclusion $\nu^*(H)=\tau^*(H)$ for any hypergraph follows from the linear programming (LP)-duality.
For $n$-vertex $k$-graphs we trivially have $\nu^*(H)=\tau^*(H)\le \frac{n}{k}$.

We construct another $(1,b)$-graph $\tilde{H}$ on $[\frac{nf}{b}]\cup V$ with edge set $E(\tilde{H})=\{\{i\}\cup e:e\in H_i$ for all $i\in[nf/b]\}$.
A $(1,k-1)$-subset $S$ of $V(\tilde{H})$ contains one vertex in $[\frac{nf}{b}]$ and $k-1$ vertices in $V$.
Let $\delta_{1,k-1}(\tilde{H}):=\min\{\deg_{\tilde{H}}(S)$: $S$ is a $(1,k-1)$-subset of $V(\tilde{H})$$\}$ where $\deg_{\tilde{H}}(S)$ denotes the number of edges in $\tilde{H}$ containing $S$.

The proof of the following claim is by now a standard argument on fractional matchings and covers.

\begin{claim}\label{frac}
If each $H_{F_i}'$ contains a perfect fractional matching for $i\in[\frac{n}{b}]$, then the auxiliary $(1,b)$-graph $H_{\textbf{G}}$ of $\textbf{G}$ contains a perfect fractional matching.
\end{claim}

\begin{proof}
By the duality theorem, we transform the maximum fractional matching problem into the minimum fractional cover problem.
Since $\tau^*(H_{\emph{\textbf{G}}})=\nu^*(H_{\emph{\textbf{G}}})\leq\frac{n}{b}$, it suffices to show that $\tau^*(H_{\emph{\textbf{G}}})\geq\frac{n}{b}$ to obtain $\nu^*(H_{\emph{\textbf{G}}})=\frac{n}{b}$.
Let $\omega$ be the minimum fractional cover of $H_{\emph{\textbf{G}}}$ and take $i_1 \in [n/b]$ such that $\omega(i_1):=\min_{i\in [n/b]}\omega(i)$.
We may assume that $\omega(i_1)=1-x<1$, since otherwise $\omega([n/b])\geq\frac{n}{b}$ and we are done.
By definition we get $\omega(e)\geq1-\omega(i_1)= x$ for every $e\in H'_{F_{i_1}}$.
We define a new weight function $\omega'$ on $V$ by setting $\omega'(v)=\frac{\omega(v)}{x}$ for every vertex $v\in V$.
Thus, $\omega'$ is a fractional cover of $H_{F_{i_1}}'$ because for each $e\in H'_{F_{i_1}}$, $\omega'(e)=\frac{\omega(e)}{x}\geq1$.
Recall that $H'_{F_{i_1}}$ has a perfect fractional matching, and thus $\omega'(V)\geq\tau^*(H_{F_{i_1}}')\geq\frac{n}{b}$ which implies that $\omega(V)\geq\frac{xn}{b}$. Therefore,

\begin{center}
$\omega([\frac{n}{b}]\cup V)\geq(1-x)\frac{n}{b}+\frac{xn}{b}=\frac{n}{b}$.
\end{center}
Hence, $\tau^*(H_{\emph{\textbf{G}}})=\frac{n}{b}$, i.e. $H_{\emph{\textbf{G}}}$ contains a perfect fractional matching.
\end{proof}



In this section, given an $n$-vertex D$k$-graph system $\emph{\textbf{G}}$, we shall construct an auxiliary $(1,b)$-graph $H_{\emph{\textbf{G}}}$ of $\emph{\textbf{G}}$ and a sequence of random subgraphs of $H_{\emph{\textbf{G}}}$.
Then, we use the properties of them to get a ``near regular'' spanning subgraph for the sake of applying Lemma \ref{pipp}.

The proof is based on a two-round randomization which is already used in \cite{MR2915641,2020A,2020B}.
Since we work with balanced $(1,b)$-graphs, we need to make sure that each random graph is balanced.
In order to achieve this we modify the randomization process by fixing an arbitrarily small and balanced set $S\subseteq V(H_{\emph{\textbf{G}}})$.
This is done in Fact \ref{balanced}.

Let $H_{\emph{\textbf{G}}}$ be the auxiliary $(1,b)$-graph of $\emph{\textbf{G}}$ with partition classes $A$, $B$ and $b|A|=|B|$ where $A$ is the color set and $B=V$.
Let $S\subseteq V(H_{\emph{\textbf{G}}})$ be a set of vertices such that $|S\cap A|=n^{0.99}/b$ and $|S\cap B|=n^{0.99}$.
The desired subgraph $H''$ is obtained by two rounds of randomization.
As a preparation to the first round, we choose every vertex randomly and uniformly with probability $p=n^{-0.9}$ to get a random subset $R$ of $V(H_{\emph{\textbf{G}}})$.
Take $n^{1.1}$ independent copies of $R$ and denote them by $R_{i+}$, $i\in[n^{1.1}]$, i.e. each $R_{i+}$ is chosen in the same way as $R$ independently.
Define $R_{i-}=R_{i+}\setminus S$ for $i\in[n^{1.1}]$.
\begin{fact}\label{balanced}
Let $n,H_{\textbf{G}},A,B,S$ and $R_{i-}, R_{i+}$ be given as above. Then, with probability $1-o(1)$, there exist subgraphs $R_i, i\in[n^{1.1}]$, such that $R_{i-}\subseteq R_i \subseteq R_{i+}$ and $R_i$ is balanced.
\end{fact}

The following two lemmas together construct the desired sparse regular $k$-graph we need.

\begin{lemma}\label{first}
Given an $n$-vertex D$k$-graph system $\textbf{G}=\{G_1,\ldots,G_{nf/b}\}$ on $V$, let $H_{\textbf{G}}$ be the auxiliary $(1,b)$-graph of $\textbf{G}$.
For each $X\subseteq V(H_{\textbf{G}})$, let $Y_X^+:=|\{i:X\subseteq R_{i+}\}|$ and $Y_X:=|\{i:X\subseteq R_i\}|$.
Let $\tilde{H}$ be with vertex set $[\frac{nf}{b}]\cup V$ and edge set $E(\tilde{H})=\{\{i\}\cup e:e\in G_i$ for all $i\in[nf/b]\}$.
Then with probability at least $1-o(1)$, we have
\begin{enumerate}
  \item $|R_i|=(1/b+1+o(1))n^{0.1}$ for all $i\in[n^{1.1}]$.
  \item $Y_{\{v\}}=(1+o(1))n^{0.2}$ for $v\in V(H_{\textbf{G}})\setminus S$ and $Y_{\{v\}}\leq(1+o(1))n^{0.2}$ for $v\in S$.
  \item $Y_{\{u,v\}}\leq2$ for all $\{u,v\}\subseteq V(H_{\textbf{G}})$.
  \item $Y_e\leq1$ for all $e\in E(H_{\textbf{G}})$.
  \item Suppose that $V(R_i)=C_i\cup V_i$, we have $\delta_{1,d}(\tilde{H}[\bigcup_{j\in C_i}[(j-1)f+1,jf]\cup V_i])\geq(c_{d,F}^{\rm cov,*}+\varepsilon/4)\binom{|R_{i+}\cap B|-d}{k-d}-|R_{i+}\cap B\cap S|\binom{|R_{i+}\cap B|-d-1}{k-d-1}\geq(c_{d,F}^{\rm cov,*}+\varepsilon/8)\binom{|R_i\cap B|-d}{k-d}$.
\end{enumerate}

\end{lemma}

\begin{lemma}\label{second}
Let $n, H_{\textbf{G}}, S, R_i$, $i\in[n^{1.1}]$ be given as in Lemma \ref{first} such that each $H_{\textbf{G}}[R_i]$ is a balanced $(1,b)$-graph and has a perfect fractional matching $\omega_i$.
Then there exists a spanning subgraph $H''$ of $H^*=\cup_iH_{\textbf{G}}[R_i]$ such that
\begin{itemize}
  \item $d_{H''}(v)\leq(1+o(1))n^{0.2}$ for $v\in S$,
  \item $d_{H''}(v)=(1+o(1))n^{0.2}$ for all $v\in V(H_{\textbf{G}})\setminus S$,
  \item $\Delta_2(H'')\leq n^{0.1}$.
\end{itemize}
\end{lemma}

The proofs follow the lines as in \cite{MR2915641,2020A,2020B} and thus we put them in the appendix.

\begin{proof}[Proof of Lemma \ref{rainbow cov}]
By the definition of $c_{d,F}^{\rm cov,*}$, Lemma \ref{first} (5) and Claim \ref{frac}, there exists a perfect fractional matching $\omega_i$ in every subgraph $H_{\emph{\textbf{G}}}[R_i]$, $i\in[n^{1.1}]$.
By Lemma \ref{second}, there is a spanning subgraph $H''$ of $H^*=\cup_iH_{\emph{\textbf{G}}}[R_i]$ such that $d_{H''}(v)\leq(1+o(1))n^{0.2}$ for each $v\in S$, $d_{H''}(v)=(1+o(1))n^{0.2}$ for all $v\in V(H_{\emph{\textbf{G}}})\setminus S$ and $\Delta_2(H'')\leq n^{0.1}$.
Hence, by Lemma \ref{pipp} (by setting $D=n^{0.2}$), $H''$ contains a cover of at most $\frac{n+n/b}{1+b}(1+a)$ edges which implies that $H''$ contains a matching of size at least $\frac{n+n/b}{1+b}(1-a(1+b-1))$, where $a$ is a constant satisfying $0<a<\phi/(1+b-1)$.
Hence $H_{\emph{\textbf{G}}}$ contains a matching covering all but at most $\phi(n+n/b)$ vertices.
\end{proof}

\section{The proof of Theorem \ref{general}}

\begin{proof}
Suppose that $\frac{1}{n}\ll\phi\ll\gamma_1\ll\gamma\ll\varepsilon'\ll\varepsilon$ where $\varepsilon', \gamma$, $\gamma_1$ are defined in Lemma \ref{absorption} and $\phi$, $\varepsilon$ in Lemma \ref{rainbow cov}.
Let $H_{\emph{\textbf{G}}}$ be the auxiliary $(1,b)$-graph of $\emph{\textbf{G}}$.
By Lemma \ref{absorption}, we get a matching $M$ in $H_{\textbf{\emph{G}}}$ of size at most $2\gamma(m-1)n$ such that for every balanced set $U\subseteq [n/b]\cup V\setminus V(M)$ of size at most $\gamma_1n$, $V(M)\cup U$ spans a matching in $H_{\emph{\textbf{G}}}$.
Let $\textbf{\emph{G}}'=\{G_1',\ldots,G_{nf/b}'\}$ be the induced D$k$-graph system of $\textbf{\emph{G}}$ on $V'$ where $V':=V\setminus V(M)$.
Denote the subsystem of $\textbf{\emph{G}}'$ by $\textbf{\emph{G}}'_I=\{G_i'\mid i\in I=[nf/b]\backslash \bigcup_{j\in V(M)\cap[n/b]}[(j-1)f+1,jf]\}$.
We still have $\delta_d(G_i')\geq (\max\{c_{d,F}^{\rm abs,*},c_{d,F}^{\rm cov,*}\}+\frac{\varepsilon}{2})\binom{n-d}{k-d}$ for $i\in I$, since $2\gamma(m-1) n\binom{n-d-1}{k-d-1}\leq\frac{\varepsilon}{2}\binom{n-d}{k-d}$.
Then, we construct the new auxiliary $(1,b)$-graph $H_{\textbf{\emph{G}}_I'}$ of $\textbf{\emph{G}}_I'$.

By Lemma \ref{rainbow cov}, $H_{\textbf{\emph{G}}_I'}$ contains a matching $M_1$ covering all but at most $\phi|V'|\leq\phi (n+n/b)$ vertices.
Suppose $W_1=[n/b]\cup V\setminus(V(M)\cup V(M_1))$, hence $|W_1|\leq \phi (n+n/b)\leq \gamma_1n$ and $W_1$ is balanced.
By Lemma \ref{absorption}, $V(M)\cup W_1$ spans a matching $M_2$ in $H_{\emph{\textbf{G}}}$ and therefore $M_1\cup M_2$ is a perfect matching in $H_{\emph{\textbf{G}}}$, which yields a rainbow $F$-factor in $\textbf{\emph{G}}$.
\end{proof}

In the next few sections we prove our results in Section 1 (Theorems~\ref{direct} --~\ref{main}), and by Theorem~\ref{general}, it suffices to specify the $\delta_d^*$ we use and bound the parameters $c_{d,F}^{\rm abs,*}$ and $c_{d,F}^{\rm cov,*}$.
Note also that we will not present a proof of Theorem~\ref{3}, as it follows from either of the two directed extensions.

\section{The proofs of Theorems~\ref{direct} and~\ref{semi}}
\label{direct section}
For Theorem~\ref{direct}, to apply Theorem~\ref{general}, we set $\delta_1^*$ as the minimum out-degree $\delta^+$.
Given that $\delta^+(D_i)\geq(1-\frac{1}{k}+\varepsilon)n$ for $i\in[\frac{n}{k}\binom{k}{2}]$,
by Theorem \ref{general}, it remains to prove that $c_{1,T_k}^{\rm abs,+}, c_{1,T_k}^{\rm cov,+}\leq1-\frac{1}{k}$.

Similarly, for Theorem~\ref{semi}, we set $\delta_1^*$ as the minimum semi-degree $\delta^0$.
Given $T\in \mathcal T_k$ and $\delta^0(D_i)\geq(1-\frac{1}{k}+\varepsilon)n$ for $i\in[\frac{n}{k}\binom{k}{2}]$,
by Theorem \ref{general}, it remains to prove that $c_{1,T}^{\rm abs,0}, c_{1,T}^{\rm cov,0}\leq1-\frac{1}{k}$.

Since the proofs are similar, we only show that $c_{1,T_k}^{\rm abs,+}, c_{1,T_k}^{\rm cov,+}\leq1-\frac{1}{k}$ for Theorem~\ref{direct}.
Recall that $T_k$ is the transitive tournament on $k$ vertices.
Note that $c_{1,T_k}^{\rm abs,+}\leq1-\frac{1}{k}$ is exactly (\ref{absdi}).

We partition the $n$-vertex digraph system $\emph{\textbf{D}}$ into $[n/k]$ subsystems $\emph{\textbf{D}}_1,\ldots,\emph{\textbf{D}}_{n/k}$ where $\emph{\textbf{D}}_i=\{D_{(i-1)\binom{k}{2}+1},\ldots,D_{i\binom{k}{2}}\}$.
Define $H_{T_k, i}$ as the $k$-graph which consists of rainbow copies of $T_k$ on $\emph{\textbf{D}}_i$ with color set $[(i-1)\binom{k}{2}+1,i\binom{k}{2}]$.
We shall show that each $H_{T_k, i}$ has a perfect fractional matching.
\begin{claim}\label{pfm}
For $i\in [n/k]$, $H_{T_k, i}$ has a perfect fractional matching.
\end{claim}

A $k$-\emph{complex} is a hypergraph $J$ such that every edge of $J$ has size at most $k$, $\emptyset \in J$ and is closed under inclusion, i.e. if $e \in J$ and $e'\subseteq e$ then $e' \in J$.
We refer to the edges of size $r$ in $J$ as $r$-\emph{edges} of $J$, and write $J_r$ to denote the $r$-graph on $V(J)$ formed by these edges.
We introduce the following notion of degree in a $k$-system $J$.
For any edge $e$ of $J$, the \emph{degree} $d(e)$ of $e$ is the number of ($|e|$+1)-edges $e'$ of $J$ which contains $e$ as a subset. (Note that this is not the standard notion of degree used in $k$-graphs, in which the degree of a set is the number of edges containing it.)
The \emph{minimum} $r$-\emph{degree} of $J$, denoted by $\delta_r(J)$, is the minimum of $d(e)$ taken over all $r$-edges $e\in J$.
Trivially, $\delta_0(J)=|V(J)|$.
So every $r$-edge of $J$ is contained in at least $\delta_r(J)$ ($r+1$)-edges of $J$.
The \emph{degree sequence} of $J$ is
\[
\delta(J)=(\delta_0(J),\delta_1(J),\ldots,\delta_{k-1}(J)).
\]

 \begin{lemma}[Lemma 3.6, \cite{MR3290271}]\label{Keevash}
 If the complex $J$ satisfies $\delta(J)\geq (n, \frac{k-1}{k}n, \frac{k-2}{k}n,\ldots, \frac{1}{k}n)$, then $J_k$ contains a perfect fractional matching.
 \end{lemma}

To prove Claim~\ref{pfm} we construct the \emph{clique k-complex} $J^i$ for each $\textbf{\emph{D}}_i$, $i\in[n/k]$, which has vertex set $V$ and edge set $E(J^i_r)$ where each edge is a rainbow $T_r$  with color set $[(i-1)\binom{k}{2}+1, (i-1)\binom{k}{2}+\binom{r}{2}]$ for each $r\in[k]$ and $i\in[n/k]$ .
Note that for each $i$, the top level $J_k^i$ is exactly $H_{T_k, i}$.
By the out-degree condition, we get
\[
\delta(J^i)\geq(n,(1-{1}/{k}+\varepsilon)n,\ldots,({1}/{k}+\varepsilon)n).
\]
Therefore Claim~\ref{pfm} follows from Lemma~\ref{Keevash} and we are done.

For completeness, we include the short proof of Lemma~\ref{Keevash} given in \cite{MR3290271}.
Given points $\mathbf{x}_1,\ldots,\mathbf{x}_s\in \mathbb{R}^d$, we define their \emph{positive cone} as
$PC(\textbf{x}_1,\ldots,\textbf{x}_s):=\{\sum_{j\in[s]}\lambda_j\textbf{x}_j:\lambda_1,\ldots,\lambda_s\geq0\}$.
Recall that $V$ is the $n$-vertex set, for any $S\subseteq V$, the characteristic vector $\chi(S)$ of $S$ is the binary vector in $\mathbb{R}^n$ such that $\chi(S)_i=1$ if and only if $i\in S$.
Given a $k$-graph $H$, if $H$ has a perfect fractional matching $\omega$, then $\textbf{1}\in PC(\chi(e):e\in H)$, since $\sum_{e\in H}\omega(e)\chi(e)=\textbf{1}$.
The well-known Farkas' Lemma reads as follows.
\begin{lemma}[Farkas' Lemma]\label{farkas}
Suppose $\textbf{v}\in \mathbb{R}^n\backslash PC(Y)$ for some finite set $Y\subseteq \mathbb{R}^n$. Then there is some $\textbf{a}\in \mathbb{R}^n$ such that $\textbf{a}\cdot\textbf{y}\geq0$ for every $\textbf{y}\in Y$ and $\textbf{a}\cdot\textbf{v}<0$.
\end{lemma}

\begin{proof}[Proof of Lemma \ref{Keevash}]
Suppose that $J_k$ does not contain a perfect fractional matching, this means that $\textbf{1}\notin PC(\chi(e):e\in J_k)$.
Then by Lemma \ref{farkas}, there is some $\textbf{a}\in \mathbb{R}^n$ such that $\textbf{a}\cdot\textbf{1}<0$ and $\textbf{a}\cdot \chi(e)\geq0$ for every $e\in J_k$. Let $V=\{v_1,\ldots,v_n\}$ and $\textbf{a}=(a_1,\ldots,a_n)$ satisfying $a_1\leq a_2\leq\cdots\leq a_n$.

We first build a $k$-edge $e=\{v_{d_1},\ldots,v_{d_k}\}\in J_k$ such that $d_j\le \frac{j-1}kn+1, j\in[k]$ as follows.
Choose $d_1=1$ and having chosen $d_1,\dots, d_j$, the choice of $d_{j+1}$ is guaranteed by $\delta_j(J)\ge (1-j/k)n$.
As $e\in J_k$, we have $\textbf{a}\cdot \chi(e)\ge 0$.
Consider $\{S_i=\{v_i,v_{i+n/k},\ldots,v_{i+(k-1)n/k}\}:i\in[n/k]\}$ which form a partition of $V$, thus $\sum_{i\in[n/k]}\textbf{a}\cdot\chi(S_i)=\textbf{a}\cdot\textbf{1}<0$.
However, as the indices of vertices of $e$ precede those of $S_i$ one-by-one and $a_1\leq a_2\leq\cdots\leq a_n$, we have for each $i$, $\textbf{a}\cdot\chi(S_i)\ge \textbf{a}\cdot\chi(e)\ge 0$, a contradiction.
\end{proof}


\section{The proof of Theorem \ref{rpartite}}
For Theorem~\ref{rpartite}, we set $\delta_1^*$ as the minimum partite-degree $\delta'$.
Let $\emph{\textbf{G}}=\{G_1,\ldots,G_{n\binom{k}{2}}\}$ be a collection of $k$-partite graphs with a common partition $V_1,\ldots,V_k$ each of size $n$ such that $\delta'(G_i)\geq(1-1/k+\varepsilon)n$ for $i\in[n\binom{k}{2}]$.
By Theorem \ref{general}, it remains to prove that $c_{1, K_k}^{\rm abs,'}, c_{1,K_k}^{\rm cov,'}\leq1-\frac{1}{k}$.

The conclusion that $c_{1,K_k}^{\rm abs,'}\leq1-\frac{1}{k}$ can be similarly derived as (\ref{absdi}) in Section~\ref{Rainbow Absorption Method} and thus omitted.
We next show how to obtain $c_{1,K_k}^{\rm cov,'}\leq1-\frac{1}{k}$.

Let $H$ be a $k$-graph and let $\mathcal{P}$ be a partition of $V(H)$. Then we say a set $S\subseteq V(H)$ is $k$-\emph{partite} if it has one vertex in any part of $\mathcal{P}$ and that $H$ is $k$-\emph{partite} if every edge of $H$ is $k$-partite. Let $V$ be a set of vertices, let $\mathcal{P}$ be a partition of $V$ into $k$ parts $V_1,\ldots, V_k$ and let $J$ be a $k$-partite $k$-system on $V$.
For each $0\leq j\leq k-1$ we define the \emph{partite minimum j-degree} $\delta^*_j(J)$ as the largest $m$ such that any $j$-edge $e$ has at least $m$ extensions to a $(j+1)$-edge in any part not used by $e$, i.e.
\[\delta^*_j(J):=\min_{e\in J_j}\min_{e\cap V_i=\emptyset}|\{v\in V_i:e\cup\{v\}\in J\}|.\]
The partite degree sequence is $\delta^*(J)=(\delta^*_0(J),\ldots,\delta^*_{k-1}(J))$.

To obtain $c_{1,K_k}^{\rm cov,'}$, we use the following lemma, which is a special case of~\cite[Lemma 7.2]{MR3290271}.
Again we present the short proof of Lemma~\ref{mul} given in \cite{MR3290271}.
\begin{lemma}\label{mul}
Let $V$ be a set partitioned into $k$ parts $V_1,\ldots,V_k$ each of size $n$ and let $J$ be a $k$-partite $k$-system on $V$ such that
\[\delta^*(J)\geq \left(n, \frac{(k-1)n}{k},\frac{(k-2)n}{k},\ldots,\frac{n}{k} \right).\]
Then $J_k$ contains a perfect fractional matching.
\end{lemma}

\begin{proof}
Suppose that $J_k$ does not contain a perfect fractional matching, this means that $\textbf{1}\notin PC(\chi(e):e\in J_k)$.
Then by Lemma \ref{farkas}, there is some $\textbf{a}\in \mathbb{R}^{kn}$ such that $\textbf{a}\cdot\textbf{1}<0$ and $\textbf{a}\chi(e)\geq0$ for every $e\in J_k$. Let $V=\{v_{1,1},\ldots,v_{1,n},v_{2,1},\ldots,v_{2,n},\ldots,v_{k,1},\ldots,v_{k,n}\}$ and
\[
\textbf{a}=(a_{1,1},\ldots,a_{1,n},a_{2,1}, \ldots, a_{2,n}, \ldots, a_{k,1},\ldots,a_{k,n})\]
satisfying that $a_{j,\frac{(j-1)n}{k}+1}\leq\cdots\leq a_{j,n}\leq a_{j,1}\leq\ldots\leq a_{j,\frac{(j-1)n}{k}}$ for each $j\in[k]$.

We first build a $k$-edge $e=\{v_{1,d_1},v_{2,d_2}\ldots,v_{k,d_k}\}\in J_k$ such that $d_j\le \frac{j-1}kn+1, j\in[k]$ as follows.
Choose $d_1=1$ and having chosen $d_1,\dots, d_j$, the choice of $d_{j+1}$ is guaranteed by $\delta_j(J)\ge (1-j/k)n$.
As $e\in J_k$, we have $\textbf{a}\cdot \chi(e)\ge 0$.
Consider $\{S_i=\{v_{1,i},v_{2,i},\ldots,v_{k,i}\}:i\in[n]\}$ which forms a partition of $V$, thus $\sum_{i\in[n]}\textbf{a}\cdot\chi(S_i)=\textbf{a}\cdot\textbf{1}<0$.
However, we have for each $i$, $\textbf{a}\cdot\chi(S_i)\ge \textbf{a}\cdot\chi(e)\ge 0$, as $a_{j,\frac{(j-1)n}{k}+1}\leq\cdots\leq a_{j,n}\leq a_{j,1}\leq\ldots\leq a_{j,\frac{(j-1)n}{k}}$ for each $j\in[k]$, a contradiction.
\end{proof}
We partition the $kn$-vertex $k$-partite graph system $\emph{\textbf{G}}$ on $V$ into $n$ subsystems $\emph{\textbf{G}}_1,\ldots,\emph{\textbf{G}}_{n}$ where $\emph{\textbf{G}}_i=\{G_{(i-1)\binom{k}{2}+1},\ldots,G_{i\binom{k}{2}}\}$ for $i\in[n]$.
The \emph{clique k-complex} $J^i$ of a $k$-partite graph system $\textbf{\emph{G}}_i$ is with vertex set $V$ and edge set $E(J^i_r)$ where each edge is a rainbow $K_r$ with color set $[(i-1)\binom{k}{2}+1,(i-1)\binom{k}{2}+\binom{r}{2}]$ for each $r\in[k]$ and $i\in[n]$.
In $\emph{\textbf{G}}$, by the degree condition, we get for each $i\in[n]$,

\[
\delta^*(J^i)\geq \left(n,\left(1-\frac{1}{k}+\varepsilon\right)n,\ldots,\left(\frac{1}{k}+\varepsilon\right)n\right).
\]

By Lemma \ref{mul}, $J_k^i$ contains a perfect fractional matching. Therefore $c_{1,K_k}^{\rm cov,'}\leq1-\frac{1}{k}$.
\section{The proof of Theorem \ref{main}}
Note that $c_{d,F}^{\rm{abs},d}\leq\frac{1}{2}$ is given in (\ref{absmatching}).
By Definition \ref{cov}, we trivially have $c_{d,F}^{\rm{cov},d}\leq c_{k,d}$ where $F$ is an edge.
By Theorem \ref{general}, the proof of Theorem \ref{main} is completed.

\section{Concluding Remarks}

In this paper we studied the rainbow version of clique-factor problems in graph and hypergraph systems.
The most desirable question is to prove an exact version of the rainbow Hajnal--Szemer\'edi theorem, which we put as a conjecture here.

\begin{conjecture}
Let $\textbf{G}=\{G_1, G_2,\ldots, G_{\frac{n}{t}\binom{t}{2}}\}$ be an $n$-vertex graph system.
If $\delta(G_i)\geq(1-\frac{1}{t})n$ for $i\in[\frac{n}{t}\binom{t}{2}]$, then $\textbf{G}$ contains a rainbow $K_t$-factor.
\end{conjecture}
Soon after this manuscript was released, Montgomery, M\"{u}yesser and Pehova \cite{montgomery2021transversal} independently proved more general results, that is, they determined asymptotically optimal minimum degree conditions for $F$-factor transversals and spanning tree transversals in graph systems.

\section{Acknowledgement}

We thank the anonymous referees for detailed feedback that improved the presentation of the paper.
 \newpage
\bibliographystyle{plain}
\bibliography{Bibte}

 \begin{appendix}

 \section{The postponed proofs}

 Below we restate and prove Fact \ref{balanced}, Lemma~\ref{first} and Lemma~\ref{second}.

\noindent \textbf{Fact 1.} Let $n,H_{\emph{\textbf{G}}},A,B,S$ and $R_{i-}, R_{i+}$ be given as above. Then, with probability $1-o(1)$, there exist subgraphs $R_i, i\in[n^{1.1}]$, such that $R_{i-}\subseteq R_i \subseteq R_{i+}$ and $R_i$ is balanced.
\begin{proof}
Recall that $|A|=n/b, |B|=n, |S\cap A|=n^{0.99}/b$ and $|S\cap B|=n^{0.99}$, thus
\[\mathbb{E}[|R_{i+}\cap A|]=n^{0.1}/b,\]
\[\mathbb{E}[|R_{i+}\cap A\cap S|]=n^{0.09}/b,\]
\[\mathbb{E}[|R_{i+}\cap B|]=n^{0.1},\]
\[\mathbb{E}[|R_{i+}\cap B\cap S|]=n^{0.09}.\]
By Lemma \ref{chernoff1}, we have
\[\mathbb{P}[||R_{i+}\cap A|-n^{0.1}/b|\geq n^{0.08}]\leq e^{-\Omega(n^{0.06})},\]
\[\mathbb{P}[||R_{i+}\cap A\cap S|-n^{0.09}/b|\geq n^{0.08}]\leq e^{-\Omega(n^{0.07})},\]
\[\mathbb{P}[||R_{i+}\cap B|-n^{0.1}|\geq n^{0.08}]\leq e^{-\Omega(n^{0.06})},\]
\[\mathbb{P}[||R_{i+}\cap B\cap S|-n^{0.09}/b|\geq n^{0.08}]\leq e^{-\Omega(n^{0.07})}.\]

Thus, with probability $1-o(1)$, for all $i\in[n^{1.1}]$,
\[|R_{i+}\cap A|\in[n^{0.1}/b-n^{0.08},n^{0.1}/b+n^{0.08}],\]
\[|R_{i+}\cap A\cap S|=(1+o(1))n^{0.09}/b,\]
\[|R_{i+}\cap B|\in[n^{0.1}-n^{0.08},n^{0.1}+n^{0.08}],\]
\[|R_{i+}\cap B\cap S|=(1+o(1))n^{0.09}.\]
Therefore, $|b|R_{i+}\cap A|-|R_{i+}\cap B||\leq (b+1)n^{0.08}<\min\{|R_{i+}\cap A\cap S|,|R_{i+}\cap B\cap S|\}$. Hence, with probability $1-o(1)$, $R_i$ can be balanced for $i\in[n^{1.1}]$.
\end{proof}
\noindent \textbf{Lemma 5.3.}
Given an $n$-vertex D$k$-graph system $\emph{\textbf{G}}=\{G_1,\ldots,G_{nf/b}\}$ on $V$, let $H_{\emph{\textbf{G}}}$ be the auxiliary $(1,b)$-graph of $\emph{\textbf{G}}$.
For each $X\subseteq V(H_{\emph{\textbf{G}}})$, let $Y_X^+:=|\{i:X\subseteq R_{i+}\}|$ and $Y_X:=|\{i:X\subseteq R_i\}|$.
Let $\tilde{H}$ be with vertex set $[\frac{nf}{b}]\cup V$ and edge set $E(\tilde{H})=\{\{i\}\cup e:e\in G_i$ for all $i\in[nf/b]\}$.
Then with probability at least $1-o(1)$, we have
\begin{enumerate}
  \item $|R_i|=(1/b+1+o(1))n^{0.1}$ for all $i\in[n^{1.1}]$.
  \item $Y_{\{v\}}=(1+o(1))n^{0.2}$ for $v\in V(H_{\emph{\textbf{G}}})\setminus S$ and $Y_{\{v\}}\leq(1+o(1))n^{0.2}$ for $v\in S$.
  \item $Y_{\{u,v\}}\leq2$ for all $\{u,v\}\subseteq V(H_{\emph{\textbf{G}}})$.
  \item $Y_e\leq1$ for all $e\in E(H_{\emph{\textbf{G}}})$.
  \item Suppose that $V(R_i)=C_i\cup V_i$, we have $\delta_{1,d}(\tilde{H}[\bigcup_{j\in C_i}[(j-1)f+1,jf]\cup V_i])\geq(c_{d,F}^{\rm cov,*}+\varepsilon/4)\binom{|R_{i+}\cap B|-d}{k-d}-|R_{i+}\cap B\cap S|\binom{|R_{i+}\cap B|-d-1}{k-d-1}\geq(c_{d,F}^{\rm cov,*}+\varepsilon/8)\binom{|R_i\cap B|-d}{k-d}$.
\end{enumerate}

      \begin{proof}

Note that
$\mathbb{E}[|R_{i+}|]=(1/b+1)n^{0.1}, \mathbb{E}[|R_{i-}|]=((1/b+1)n-(1/b+1)n^{0.99})n^{-0.9}=(1/b+1)n^{0.1}-(1/b+1)n^{0.09}$.
By Lemma \ref{chernoff1}, we have
\[\mathbb{P}[\mid|R_{i+}|-n^{0.1}(1/b+1)\mid\geq n^{0.095}]\leq e^{-\Omega(n^{0.09})},\]
\[\mathbb{P}[\mid|R_{i-}|-((1/b+1)n^{0.1}-(1/b+1)n^{0.09})\mid\geq n^{0.095}]\leq e^{-\Omega(n^{0.09})}.\]
Hence, with probability at least $1-O(n^{1.1})e^{-\Omega(n^{0.09})}$, for the given sequence $R_i$ in Fact \ref{balanced}, $i\in[n^{1.1}]$, satisfying $R_{i-}\subseteq R_i \subseteq R_{i+}$, we have $|R_i|=(1/b+1+o(1))n^{0.1}$.

For each $X\subseteq V(H_{\emph{\textbf{G}}})$, let $Y_X^+:=|\{i:X\subseteq R_{i+}\}|$ and $Y_X:=|\{i:X\subseteq R_i\}|$. Note that the random variables $Y_X^+$ have binomial distributions $Bi(n^{1.1},n^{-0.9|X|})$ with expectations $n^{1.1-0.9|X|}$ and $Y_X\leq Y_X^+$. In particular, for each $v\in V(H_{\emph{\textbf{G}}})$, $\mathbb{E}[Y_{\{v\}}^+]=n^{0.2}$, by Lemma \ref{chernoff1}, we have
\[\mathbb{P}[\mid|Y_{\{v\}}^+|-n^{0.2}\mid\geq n^{0.19}]\leq e^{-\Omega(n^{0.18})}.\]
Hence, with probability at least $1-O(n)e^{-\Omega(n^{0.18})}$, we have $Y_{\{v\}}=(1+o(1))n^{0.2}$ for $v\in V(H_{\emph{\textbf{G}}})\setminus S$ and $Y_{\{v\}}\leq(1+o(1))n^{0.2}$ for $v\in S$.

Let $Z_{p,q}=|X\in\binom{V(H_{\emph{\textbf{G}}})}{p}:Y_X^+\geq q|$. Then,
\[\mathbb{E}[Z_{p,q}]\leq\binom{\frac{n}{b}+n}{p}\binom{n^{1.1}}{q}(n^{-0.9pq})\leq Cn^{p+1.1q-0.9pq}.\]
Hence, by Markov's inequality we have
\[\mathbb{P}[Z_{2,3}=0]=1-\mathbb{P}[Z_{2,3}\geq1]\geq1-\mathbb{E}[Z_{2,3}]=1-o(1),\]
\[\mathbb{P}[Z_{1+b,2}=0]=1-\mathbb{P}[Z_{1+b,2}\geq1]\geq1-\mathbb{E}[Z_{1+b,2}]=1-o(1),\]
i.e. with probability at least $1-o(1)$, every pair $\{u,v\}\subseteq V(H_{\emph{\textbf{G}}})$ is contained in at most two sets $R_{i+}$, and every edge is contained in at most one set $R_{i+}$. Thus, the conclusions also hold for $R_i$.

Fix a $(1,d)$-subset $D\subseteq V(\tilde{H})$ and let $N_D(\tilde{H})$ be the neighborhood of $D$ in $\tilde{H}$.
Recall that $R$ is obtained by choosing every vertex randomly and uniformly with probability $p=n^{-0.9}$,
let $DEG_D$ be the number of edges $\{f|f\subseteq R$ and $f\in N_D(\tilde{H})\}$.
Therefore $DEG_D=\sum_{f\in N_D(\tilde{H})}X_f$, where $X_f=1$ if $f$ is in $R$ and 0 otherwise. We have
\[\mathbb{E}[DEG_D]=d_{\tilde{H}}(D)\times(n^{-0.9})^{k-d}\geq(c_{d,F}^{\rm cov,*}+\varepsilon)\binom{n-d}{k-d}n^{-0.9(k-d)}\]
\[\geq(c_{d,F}^{\rm cov,*}+\varepsilon/3)\binom{|R\cap B|-d}{k-d}=\Omega(n^{0.1(k-d)}).\]

For two distinct intersecting edges $f_i,f_j\in N_D(\tilde{H})$ with $|f_i\cap f_j|=\ell$ for $\ell\in[k-d-1]$, the probability that both of them are in $R$ is
\[\mathbb{P}[X_{f_i}=X_{f_j}=1]=p^{2(k-d)-\ell},\]
for any fixed $\ell$, we have
\[\Delta=\sum_{f_i\cap f_j\neq\emptyset}\mathbb{P}[X_{f_i}=X_{f_j}=1]\leq\sum_{\ell=1}^{k-d-1}p^{2(k-d)-\ell}\binom{n-d}{k-d}\binom{k-d}{\ell}\binom{n-k}{k-d-\ell}\]
\[\leq\sum_{\ell=1}^{k-d-1}p^{2(k-d)-\ell}O(n^{2(k-d)-\ell})=O(n^{0.1(2(k-d)-1)}).\]
Applying Lemma \ref{Janson} with $\Gamma=B$, $\Gamma_p=R\cap B$ and $M=N_{\tilde{H}}(D)$(a family of $(k-d)$-sets) , we have
\[\mathbb{P}[DEG_D\leq(1-\varepsilon/12)\mathbb{E}[DEG_D]]\leq e^{-\Omega((\mathbb{E}[DEG_D])^2/\Delta)}=e^{-\Omega(n^{0.1})}.\]
Therefore by the union bound, with probability $1-o(1)$, for all $(1,d)$-subsets $D\subseteq V(\tilde{H})$, we have
\[DEG_D>(1-\varepsilon/12)\mathbb{E}[DEG_D]\geq(c_{d,F}^{\rm cov,*}+\varepsilon/4)\binom{|R\cap B|-d}{k-d}.\]
Summarizing, with probability $1-o(1)$, for the sequence $R_i, i\in[n^{1.1}]$, satisfying $R_{i-}\subseteq R_i \subseteq R_{i+}$, all of the following hold.
\begin{enumerate}
  \item $|R_i|=(1/b+1+o(1))n^{0.1}$ for all $i\in[n^{1.1}]$.
  \item $Y_{\{v\}}=(1+o(1))n^{0.2}$ for $v\in V(H_{\emph{\textbf{G}}})\setminus S$ and $Y_{\{v\}}\leq(1+o(1))n^{0.2}$ for $v\in S$.
  \item $Y_{\{u,v\}}\leq2$ for all $\{u,v\}\subseteq V(H_{\emph{\textbf{G}}})$.
  \item $Y_e\leq1$ for all $e\in E(H_{\emph{\textbf{G}}})$.
  \item $DEG_D^{(i)}\geq(c_{d,F}^{\rm cov,*}+\varepsilon/4)\binom{|R_{i+}\cap B|-d}{k-d}$ for all $(1,d)$-subsets of $D\subseteq V(\tilde{H})$ and $i\in[n^{1.1}]$.
\end{enumerate}
Thus, by property (5) above, we conclude that suppose $V(R_i^+)=C_i^+\cup V_i^+$ and $V(R_i)=C_i\cup V_i$, the following holds.
\[\delta_{1,d}(\tilde{H}[\bigcup_{j\in C_i^+}[(j-1)f+1,jf]\cup V_i^+]])\geq(c_{d,F}^{\rm cov,*}+\varepsilon/4)\binom{|R_{i+}\cap B|-d}{k-d}.\]

After the modification, we still have
\[\delta_{1,d}(\tilde{H}[\bigcup_{j\in C_i}[(j-1)f+1,jf]\cup V_i]])\geq(c_{d,F}^{\rm cov,*}+\varepsilon/4)\binom{|R_{i+}\cap B|-d}{k-d}-|R_{i+}\cap B\cap S|\binom{|R_{i+}\cap B|-d-1}{k-d-1}\]
\[\geq(c_{d,F}^{\rm cov,*}+\varepsilon/8)\binom{|R_i\cap B|-d}{k-d}. \qedhere
\]
\end{proof}

\noindent \textbf{Lemma 5.4.}
Let $n, H_{\emph{\textbf{G}}}, S, R_i$, $i\in[n^{1.1}]$ be given as in Lemma \ref{first} such that each $H_{\emph{\textbf{G}}}[R_i]$ is a balanced $(1,b)$-graph and has a perfect fractional matching $\omega_i$.
Then there exists a spanning subgraph $H''$ of $H^*=\cup_iH_{\emph{\textbf{G}}}[R_i]$ such that
\begin{itemize}
  \item $d_{H''}(v)\leq(1+o(1))n^{0.2}$ for $v\in S$,
  \item $d_{H''}(v)=(1+o(1))n^{0.2}$ for all $v\in V(H_{\emph{\textbf{G}}})\setminus S$,
  \item $\Delta_2(H'')\leq n^{0.1}$.
\end{itemize}

      \begin{proof}
By the condition that each $H_{\emph{\textbf{G}}}[R_i]$ has a perfect fractional matching $\omega_i$, we select a generalized binomial subgraph $H''$ of $H^*$ by independently choosing each edge $e$ with probability $\omega_{i_e}(e)$ where $i_e$ is the index $i$ such that $e\in H_{\emph{\textbf{G}}}[R_i]$.
Recall that property (4) guarantees the uniqueness of $i_e$.

For $v\in V(H'')$, let $I_v=\{i:v\in R_i\}$, $E_v=\{e\in H^*:v\in e\}$ and $E_v^i=E_v\cap H_{\emph{\textbf{G}}}[R_i]$, then $E_v^i, i\in I_v$ forms a partition of $E_v$ and $|I_v|=Y_{\{v\}}$.
Hence, for $v\in V(H'')$,
\[d_{H''}(v)=\sum_{e\in E_v}1=\sum_{i\in I_v}\sum_{e\in E_v^i}X_e,\]
where $X_e$ is the Bernoulli random variable with $X_e=1$ if $e\in E(H'')$ and $X_e=0$ otherwise. Thus its expectation is $\omega_{i_e}(e)$.
Therefore
\[\mathbb{E}[d_{H''}(v)]=\sum_{i\in I_v}\sum_{e\in E_v^i}\omega_{i_e}(e)=\sum_{i\in I_v}1=Y_{\{v\}}.\]
Hence, $\mathbb{E}[d_{H''}(v)]=(1+o(1))n^{0.2}$ for $v\in V(H_{\emph{\textbf{G}}})\setminus S$ and $\mathbb{E}[d_{H''}(v)]\leq(1+o(1))n^{0.2}$ for $v\in S$. Now by Chernoff's inequality, for $v\in V(H_{\emph{\textbf{G}}})\setminus S$,
\[\mathbb{P}[|d_{H''}(v)-n^{0.2}|\geq n^{0.15}]\leq e^{-\Omega(n^{0.1})},\]
and for $v\in S$,
\[\mathbb{P}[d_{H''}(v)-n^{0.2}\geq n^{0.15}]\leq e^{-\Omega(n^{0.1})},\]

Taking a union bound over all vertices, we conclude that with probability $1-o(1)$, $d_{H''}(v)=(1+o(1))n^{0.2}$ for all $v\in V(H_{\emph{\textbf{G}}})\setminus S$ and $d_{H''}(v)\leq(1+o(1))n^{0.2}$ for $v\in S$.

Next, note that for distinct $u,v\in V(H_{\emph{\textbf{G}}})$,
\[d_{H''}(\{u,v\})=\sum_{e\in E_u\cap E_v}1=\sum_{i\in I_u\cap I_v}\sum_{e\in E_u^i\cap E_v^i}X_e,\]
and
\[\mathbb{E}[d_{H''}(\{u,v\})]=\sum_{i\in I_u\cap I_v}\sum_{e\in E_u^i\cap E_v^i}\omega_i(e)\leq|I_u\cap I_v|\leq2.\]
Thus, by Lemma \ref{Chernoff2},
\[\mathbb{P}[d_{H''}(\{u,v\})\geq n^{0.1}]\leq e^{-n^{0.1}},\]
then by a union bound we have $\Delta_2(H'')\leq n^{0.1}$ with probability $1-o(1)$.
\end{proof}
  \end{appendix}

\end{document}